# An introduction to functional analysis for science and engineering

David A. B. Miller

*Stanford University*

This article is a tutorial introduction to the functional analysis mathematics needed in many physical problems, such as in waves in continuous media. This mathematics takes us beyond that of finite matrices, allowing us to work meaningfully with infinite sets of continuous functions. It resolves important issues, such as whether, why and how we can practically reduce such problems to finite matrix approximations. This branch of mathematics is well developed and its results are widely used in many areas of science and engineering. It is, however, difficult to find a readable introduction that both is deep enough to give a solid and reliable grounding but yet is efficient and comprehensible. To keep this introduction accessible and compact, I have selected only the topics necessary for the most important results, but the argument is mathematically complete and self-contained. Specifically, the article starts from elementary ideas of sets and sequences of real numbers. It then develops spaces of vectors or functions, introducing the concepts of norms and metrics that allow us to consider the idea of convergence of vectors and of functions. Adding the inner product, it introduces Hilbert spaces, and proceeds to the key forms of operators that map vectors or functions to other vectors or functions in the same or a different Hilbert space. This leads to the central concept of compact operators, which allows us to resolve many difficulties of working with infinite sets of vectors or functions. We then introduce Hilbert-Schmidt operators, which are compact operators encountered extensively in physical problems, such as those involving waves. Finally, it introduces the eigenvectors or eigenfunctions for major classes of operators, and their powerful properties, and ends with the important topic of singular-value decomposition of operators. This article is written in a style that is complementary to that of standard mathematical treatments; by relegating longer proofs to a separate section, I have attempted to retain a clear narrative flow and motivation in developing the mathematical structure. Hopefully, the result is useful to a broader readership who need to understand this mathematics, especially in physical science and engineering.

## 1 Introduction

Physical scientists and engineers are typically well educated in many branches of mathematics. Sets, the various kinds of numbers, calculus, differential equations, and linear algebra (especially with finite matrices) form a typical grounding. It is not uncommon in these disciplines to encounter results from another field of mathematics when we have to work with sets of functions; this is routine in quantum mechanics, for example, which is mathematically built around the general linear algebra of operators and sets of eigenfunctions. But that field of mathematics is not itself part of the typical course sequence for such scientists and engineers. When we need to understand those results more deeply, we therefore have a problem. Recently, in understanding problems with waves [1], for example, such as meaningfully counting the number of usable communications channels between sources and receivers, that lack of understanding has led to substantial confusion and even error[1]. This "missing" field of mathematics is functional analysis.

Functional analysis is a highly developed field that is well-known to mathematicians. But, possibly because it is not generally taught to others, its literature is resolutely mathematical, erecting a higher barrier of incomprehensibility for other mortals. Its texts, like many in mathematics, tend to be dry to the point of total dessication. That may suit those who like mathematics for its own sake, and may even be regarded as a virtue in that discipline. But, for others whose motivation is to understand the mathematics they need for specific actual problems in the physical world, the lack of any narrative structure (so – where are we going with this?) and of any sense of purpose in this next definition, lemma or proof (so – do I need this, and if so why?) can make the field practically inaccessible to many others who could understand it well and use it productively.

---

[1] Indeed, the impetus for writing this review was precisely to give the necessary background for such a deeper analysis of such wave problems [1]

This article is a modest attempt to introduce at least an important and useful part of this field of functional analysis. By a careful selection of topics, by avoiding the temptation of every incidental lemma, and by relegating major proofs to the end, I hope to construct a narrative that leads us willingly through to understanding. A few of those proofs do indeed require some deep intellectual effort, especially if one is not used to certain kinds of mathematical arguments. I have, however, relegated all of these proofs to the end, so the reader can more easily follow the sequence of useful ideas. With that, this field, or at least the part I will cover, is, I believe, relatively comprehensible, and even straightforward.

Broadly speaking, functional analysis takes the kinds of results that are simple and even obvious for concepts such as the convergence of sequences of real numbers, and extends them to other situations; we can then examine especially the convergence of sequences of different vectors (possibly with large or even infinite numbers of elements) or continuous functions, and even sequences of matrices or linear operators. This extension is possible because we only need to generalize a few of the properties we find with real numbers; then we can establish the necessary techniques for convergence of these more complex vectors, functions, matrices, and operators. We can then build on those techniques to generate powerful results.

## 2   Ideas and approach

Before starting on the mathematics proper, we can summarize where we are going, why we are going there, and how we going to get there. We can do this first by looking at the key ideas, and then give a "roadmap" for this article.

### 2.1   Key ideas

We give a brief and informal introduction here the kinds of ideas and mathematical concepts we will use. We will develop all the formal definitions and resulting mathematics later.

#### 2.1.1   Vectors and functions

In this mathematics, instead of just working with individual numbers, we are going to have to work at least with "vectors". Here these vectors are almost always not the geometrical vectors with $x$, $y$, and $z$ components; instead, we could regard them as representing the values of a function at every point of interest, and we can think of them as columns of those numbers – the values of the function at each such point – which are vectors in the sense of matrix-vector algebra. Of course, immediately we see that, to represent "smooth" functions, these vectors are likely to be infinitely long. Since these vectors are representing functions, we can (and do) refer to them as functions, and we will use the terms "vectors" and "functions" interchangeably in what follows here.

#### 2.1.2   Operators

In normal matrix algebra with finite matrices, we are used to the idea that a matrix can "operate on" (that is, multiply) a vector (usually considered to be on the right of the matrix) to generate a new vector. A matrix operating in this way is technically a "linear operator". In the mathematics we develop here, we generalize this idea of a linear operator; the reader can, however, be confident that this operator can always in the end be thought of as a matrix, possibly of infinite dimensions. Importantly, with the same functional analysis mathematics, this operator could be a linear integral operator (as in a Fourier transform operator, or a Green's function, for example), kinds of operators that are common in working with continuous functions.

For this mathematics, we have to define some additional concepts for such operators, especially whether they are what is known as "compact". This "compactness" property is essentially what allows us to approximate infinite systems with appropriate finite ones. Such compactness is trivial for finite matrices (all finite matrices have this property), but other operators may not be compact.

A particularly important class of operators – "Hilbert-Schmidt" operators – comes up frequently in physical problems, and we will see these are all compact. They have another important property that allows us to define "sum rules" for various physical problems. A further important characteristic of operators is whether they are "Hermitian". For a matrix, this simply means that, if we transpose the matrix and then take the complex conjugate of all the elements, we recover the original matrix. For more general operators, we can give a correspondingly more general definition. This "Hermiticity" appears often in our applications, and yields further very useful mathematical results.

Operators can lead to very specific functions – the "eigenfunctions" – which are those functions that, when we operate on them with the operator, lead to the same function at the output, just multiplied by some "eigenvalue". This is analogous to the similar behavior for eigenvectors and eigenvalues for matrices, but these functions can have many further useful properties in our case.

### 2.1.3   Norms, metrics, and inner products

With ordinary numbers, the idea of the "size" of a number is obvious – it is just the magnitude of the number. For functions, vectors, and matrices, we need a generalization of this idea, which we call the "norm". With ordinary numbers, the "distance" between two numbers is also obvious – we subtract one from the other, and take the magnitude of that. For functions, vectors, and matrices, we need an analogous concept, which we call the "metric". We we can visualize this loosely as the length of the "vector" that joins the tips of the two "vectors" we are comparing. With that metric, we can talk about the convergence of sequences of "vectors" or sequences of functions. Both the norm and the metric can be simply defined based on another idea – the "inner product". The inner product between two vectors or functions can be viewed as a generalization of the geometrical vector idea of a vector "dot" product. The inner product is mathematically of central importance here, and often has substantial physical meaning.

### 2.1.4   Hilbert spaces

Just as we can think of ordinary geometric vectors as existing in ordinary three-dimensional "geometrical" space, we can consider our more general vectors as existing in a space (which generally may have an infinite number of dimensions). If we give that space the property of having an inner product, then (with some minor additional properties) we can call this space a "Hilbert space", which then becomes our analog of geometrical space.

### 2.1.5   Basis sets and eigenfunctions

We will see that, just as we can define unit vectors in the $x$, $y$, and $z$ directions in ordinary geometrical space and use them to describe any vector, we can similarly define "basis" vectors in our Hilbert spaces. Importantly, the eigenfunctions of many of the operators we will work with will be "complete" basis sets, which have many useful mathematical properties. A final step in our mathematical process is to establish the operation of "singular value decomposition" of an operator, which can define pairs of eigenfunctions in different Hilbert spaces; this is particularly useful, for example, for looking at problems with waves, where we may be generating waves in one "volume" from sources in another.

## 2.2   Structure of this article

With these various ideas, we can set up powerful mathematics that retains much of our understanding for finite matrices and vectors and extends it for continuous functions. Below, then, we develop this mathematics, including all the definitions of terms and the necessary theorems, and including proofs. Note that when we are defining some term below, we will write the term in *italic* font. I have including an index of all of these definitions at the end (section 13) so the reader can refer back to those as the reader progresses further into the development here.

This article is customized in its order and in the specific material included, and the overall argument presented here is original in that sense. The underlying mathematical material is, however, standard; nearly all the mathematical definitions, theorems and proofs below are found in texts in functional analysis. In particular, I have made extensive use of discussions in Kreyszig [2], Hunter and Nachtergaele [3], and Hanson and Yakovlev [4], as well as general mathematical reference material[2]. If some proof closely follows another treatment, I give that specific source. I have tried to use a consistent notation within this article, though that notation is necessarily not quite the same as any of these sources.

The construction of the argument here is also original in that it avoids introducing much mathematics that we do not need; by contrast, the standard treatments in functional analysis texts typically insist on building up through a progressive and very general set of concepts, many of which are irrelevant here and that make it more difficult to get to the ideas we do need[3]. The only actual mathematical innovations in this article are in some minor but novel notations we introduce that make the mathematics easier to use in some situations[4].

The way that I present the material here is also different in structure to the approach of typical mathematics texts. I emphasize the key points and structure of the overall argument rather than the detailed proofs. I include all the necessary proofs, but the more involved proofs are included at the end so that those details do not interrupt the overall logical and narrative flow. If this mathematics is new to you, I recommend reading through the argument here up to section 11, and returning to that section later for deeper study of the various proofs. Overall this article is by far the shortest path I know of to understanding this material. Hopefully, too, this article may help the reader follow these more comprehensive texts [2][3][4] if the reader needs yet more depth or an alternate treatment.

In this article, in section 3, I will introduce the necessary concepts from the analysis of convergence with ordinary numbers. Section 4 extends these convergence ideas to functions and vectors. Section 5 introduces Hilbert spaces, and section 6 continues by introducing key ideas for operators. Section 7 treats the eigenfunctions and eigenvalues of the most important class of operators for our purposes (compact Hermitian operators). Section 8 expands the concept of inner products, allowing ones with additional specific physical meanings, and section 9 gives the extension of this algebra to singular-value decomposition. After concluding remarks in section 10, various proofs are given in detail in section 11. After references in section 12, an index of definitions is given in section 13.

# 3  Convergence of sequences of real numbers

To think about how one vector, function, matrix or operator matrix converges on another such entity, we need to start by thinking about how sequences of numbers converge, and we generally do this by thinking about convergence of sequences of real numbers[5]. To do so, we should be clear about what a sequence of real numbers is. We should also formally introduce the ideas of norms and metrics, which are very straightforward for real numbers. (These ideas can be applied later to other entities, such as vectors, functions, matrices and operators.) Then we can formally define convergence of a sequence of real numbers and give some other important definitions. Following that, we can introduce the Bolzano-Weierstrass theorem (proved in section 11), which is a useful mathematical tool for later proofs, and is an essential concept for understanding all the ideas of convergence.

---

[2] Kreyszig [2] is a classic introductory text that is more comprehensive and readable than many texts in the field of functional analysis, but it omits explicitly discussion of Hilbert-Schmidt operators (though it does cover much of the associated mathematics). Hunter and Nachtergaele [3] is an example of a more modern and relatively complete text. Hanson and Yakovlev [4] is not so complete on the mathematics of functional analysis in itself (though it refers to additional proofs in other sources), but includes substantial discussion of applications in electromagnetism.

[3] As a result, my approach here takes up only about 10% of the corresponding length of such mathematics texts in getting to the same point. Of course, there is much other good mathematics in those texts, but that other 90% makes it much harder to understand just the mathematics we do need.

[4] In part because we may be working with operators that map from one Hilbert space to another, and those may have different inner products, we introduce the explicit notation of an "underlying inner product". We also expand the use of the Dirac notation. This is a common linear algebra notation in some fields (especially quantum mechanics), but is not common in mathematics texts. We are also able to make full use of it, through what we call an "algebraic shift", which is essentially a shift from algebra using general inner products, whose algebra is more restricted, to one that is a simple analog of matrix-algebra. Dirac notation can be regarded as just a convenient general notation for the algebra of complex matrices.

## 3.1    Sets and sequences of real numbers

Before defining a sequence, we can define a *set* of real numbers. A set can be defined just by writing out all its members. A set itself need have no properties other than having members, though we can give it additional properties. Conventionally, we write a set out by enclosing a list of its elements within curly brackets. So $\{1.7, 3.6, 2\}$ is a set that contains the three real numbers 1.7, 3.6, and 2. The order of the elements of the set (here the 3 real numbers) does not matter if we are just listing the elements, so $\{1.7, 2, 3.6\}$ means the same as $\{1.7, 3.6, 2\}$.

We often do care about the order of numbers, however. 1.7, 2, and 3.6 might be the values of some function at successive points, or they might be the $x$, $y$, and $z$ coordinates of some point in ordinary geometric space. If we care about the order, then we use a *sequence* of elements, which conventionally we write by enclosing the elements in ordinary braces, as in $(1.7, 3.6, 2)$; this is a sequence of the three real numbers 1.7, 3.6, and 2 in this specific order. So, the sequence $(1.7, 3.6, 2)$ is different from the sequence $(1.7, 2, 3.6)$. We could also write the sequence $(1.7, 3.6, 2)$ using a notation $(1.7, 3.6, 2) \equiv (x_j)$ where $x_1 = 1.7$, $x_2 = 3.6$ and $x_3 = 2$, and we presume $j$ successively takes on the values 1, 2, and 3 – i.e., $j = 1, 2, 3$.

Sequences of elements can obviously be finite in length, as in the above simple examples. However, it is common in mathematics texts on functional analysis to presume implicitly that sequences written in the form $(x_j)$ are infinite in length unless otherwise stated[6], with the variable $j$ typically then ranging over all the positive integers, or, equivalently, the "natural numbers"[7], which explicitly start at 1, excluding zero.

Below, we will also need the idea of a "subsequence". A *subsequence* is formed from a sequence by choosing some of the elements of the sequence while keeping the order of those elements the same. So $(1.7, 2)$ is a subsequence of the sequence $(1.7, 3.6, 2)$, whereas $(2, 1.7)$ is not. Note that a subsequence is also a sequence in its own right, being a set of elements in a specific order.

## 3.2    Sets and spaces

The set of (all) real numbers[8] is a set that also possesses various axiomatic properties. These are, essentially, the ordinary arithmetic properties of addition, subtraction, multiplication and division. With these additional properties, the set of real numbers is also a (mathematical) *field*, as is the set of complex numbers.

When we are operating with vectors, functions, operators or matrices, the corresponding sets may also have some axiomatic properties, especially the ideas of norms and metrics below; these sets might not have all the same algebraic properties[9] as fields, so we need to take some care in defining the properties they do have. An *(abstract) space* is a set with some additional axiomatic properties, and this is what we mean below by the term "space" in a mathematical context. A (mathematical) field is also a specific example of a space, so we can use the term "space" to cover sets with many different kinds of added attributes.

---

[5] The arguments here also work without change for complex numbers.

[6] This default assumption of infinite length is not always clear or explicit in texts on functional analysis, which can cause significant confusion in reading them.

[7] The set of all natural numbers is usually denoted by the symbol $\mathbb{N}$.

[8] The set of all real numbers is usually denoted by the symbol $\mathbb{R}$.

[9] For example, multiplication of operators or matrices is not generally commutative, and division of one vector, function or matrix by another may not have any meaning.

## 3.3    Norms and metrics

To talk about the "size" or "magnitude" of elements, such as real numbers (or, later, functions, vectors, operators and matrices), we can introduce the idea of a norm. To be useful in the mathematics that follows, a *norm* of an element $x$ in a set or space, written as $\|x\|$, is a real number, and for arbitrary elements $x$ and $y$ in the set or space, it must possess the following four additional properties[10]:

(N1)      $\|x\| \geq 0$

(N2)      $\|x\| = 0$ if and only if $x = 0$

(N3)      $\|ax\| = |a| \|x\|$ where $a$ is any *scalar* (i.e., real or complex number)

(N4)      $\|x + y\| \leq \|x\| + \|y\|$ (the *triangle inequality* for norms)

$$(1)$$

For the set or space of real numbers, we choose the norm to be simply the modulus – i.e., the norm of the real number $x$ is $|x|$ (and we can make the same choice for the norm of the set or space of complex numbers if needed). A set or space on which we have defined a norm is called a *normed space*.

A *metric* $d$ in a set or space expresses a notion of "distance" between two elements $x$ and $y$ of the set or space, and has to have four properties for arbitrary elements $x$ and $y$:

(M1)      $d$ is real-valued, finite and nonnegative (i.e., $d \geq 0$)

(M2)      $d(x, y) = 0$ if and only if $x = y$ (i.e., if and only if $x$ and $y$ are the same element)

(M3)      $d(x, y) = d(y, x)$ (symmetry)

(M4)      $d(x, y) \leq d(x, z) + d(z, y)$ (the *triangle inequality*)

$$(2)$$

For real numbers, we can choose a metric

$$d_{\mathbb{R}}(x, y) \equiv |x - y| \tag{3}$$

i.e., the modulus of the difference between the two numbers. A set or space on which we have defined a metric is called a *metric space*. A metric like this one $d_{\mathbb{R}}(x, y)$, which obviously follows directly from the definition of the norm, is sometimes called the *metric induced by the norm*.

## 3.4    Convergence of a sequence

With the ideas of metrics (and norms), the notion of convergence of a sequence of real numbers is straightforward. Generally, an (infinitely long) sequence $(x_n)$ in a metric space with metric $d$ is said to *converge* if there is an element $x$ in the space such that

$$\lim_{n \to \infty} d(x_n, x) = 0 \tag{4}$$

---

[10] Note that in the expression $\|x\| = 0$ here, the "0" is the number zero, but in the general case where $x$ may be a vector rather than a number, the expression $x = 0$ has to be taken to mean that the "0" refers to the zero vector, a vector of zero "length" or norm. This ambiguity of notation is unfortunately common in this mathematics.

where "$\lim_{n\to\infty}$" is a short-hand for "the limit as $n$ tends to infinity".

Defined in this way, this idea of convergence can later be applied to other kinds of spaces in which the elements may not be real numbers but for which we have defined some metric – specifically, the elements might later be functions, vectors, operators, or matrices. We then call $x$ the *limit* of $(x_n)$ and we can write either $\lim_{n\to\infty} x_n = x$ or the notation $x_n \to x$; both notations are equivalent to saying that $(x_n)$ converges to $x$ or has limit $x$. (If $(x_n)$ is not convergent, then by definition it is *divergent*.)

Note that $x$ must be in the set or space if it is to be a limit in this strict sense. It is quite possible for a sequence to converge to a limit that lies outside a set or space, though it would have to be "just" outside. For example, if we considered the set or space of real numbers starting from (and including) 0 and continuing up to, but not including, 1 (a set that could be written as $[0,1)$ in one common mathematical notation), the sequence $(x_n)$ with elements $x_n = 1 - 1/n$, i.e., with elements 0, ½, 2/3, ¾, 4/5, … and so on, is converging towards 1 as its limit, but 1 is not itself in the sequence. In this case, the sequence, though convergent, is not converging to a limit in the set or space (even though it is converging to a value that lies just outside the set or space). Of course, we could easily "close" this space by adding in the element 1; the resulting set of all these elements and the element 1 could then be called the *closure* of the original set. (Note that the closure is not just the additional elements required to close the set; it is all the elements of the original set plus those required to close it.)

One particularly important property of a set of elements is whether it is *bounded*. Formally, a set is bounded if, for any choices of elements $x$ and $y$ in the set, the supremum of the metric $d(x, y)$ is finite[11]. The *supremum* of a set of real numbers is the smallest number that is greater than or equal to all the elements of the set, if such an element exists. The supremum is also referred to as the least upper bound.

Often, the supremum will be the maximum element in the set, but the supremum and the maximum are not necessarily the same; a set may have a supremum that is not in the set, whereas a maximum, if it exists, would have to be in the set. For example, the infinite set of elements ½, 3/4, 7/8, 15/16, 31/32, and so on, has a supremum of 1, but the element 1 is not in this set, and it is not clear that there is a maximum element in this set; for any specific element we choose that is arbitrarily close to 1, there is another element that is even closer. We can call a sequence $(x_n)$ a *bounded sequence* if the corresponding set of points in the sequence is a bounded set.

In discussing boundedness, we sometimes also need the complementary concept of an infimum, especially if the numbers in question could be both positive and negative. The *infimum* of a set of real numbers is the largest number that is less than or equal to all the elements of the set, if such an element exists. Similarly, it is not necessarily the same as the minimum element in the set; for example, a minimum may not exist in a set, as in the set with elements ½, ¼, 1/8, 1/16, and so on, which has an infimum of 0, but may have no minimum element. (Metrics are always positive or zero, so they naturally have a lower bound of zero, so with metrics we may not need to deal with the infimum explicitly.)

## 3.5    Cauchy sequences and complete spaces

For subsequent definitions and proofs, in addition to the definition in Eq. (4), it is useful also to have a slightly different kind of notion of convergence of sequences. To set this up, we first define the idea of a *Cauchy sequence*. Specifically, in a metric space with a metric $d(x, y)$ (and the real numbers are a good example, with the metric $d_{\mathbb{R}}$ as defined in Eq. (3)), a sequence $(x_n)$ is said to be Cauchy (or to be a Cauchy sequence)[12] if for every real number $\varepsilon > 0$ (no matter how small) there is a number $N$ (a positive integer or natural number) such that,

$$\text{for every } m, n > N, \ d(x_m, x_n) < \varepsilon \tag{5}$$

---

[11] Note, incidentally, that this notion of boundedness only uses the "distance" between any two elements, not the "value" of the elements themselves. This choice is slightly less restrictive formally because it means we do not need to know the absolute size of the elements, and it is a sufficient definition for the proofs we need to construct. For the set of real numbers between 100 and 102, the supremum of this metric would be 2, not 102.

So, once we get past some point in the sequence (specifically after the $N$th element), the elements are all closer to one another than some chosen (positive) "distance" or separation $\varepsilon$; and no matter how small a distance (i.e., $\varepsilon$) we choose, there is always some $N$ (possibly different for each $\varepsilon$) for which this holds.

The distinction between a convergent sequence as in Eq. (4) and a Cauchy sequence (which by definition converges as in Eq. (5)) essentially makes no difference for us, because we can prove that

$$\text{every convergent sequence in a metric space is a Cauchy sequence} \tag{6}$$

See 11.1 "Proof (1) that every convergent sequence in a metric space is a Cauchy sequence" below.

A (metric) space is said to be *complete* if every Cauchy sequence in the space converges to a limit that is also an element of the space[13]. Of course, this is not saying that in a complete metric space every sequence is a Cauchy sequence, or even that every sequence converges, but it is saying that if every Cauchy sequence converges in the space, then the space is "complete" (by definition). A complete space therefore has to have all the limit points of Cauchy sequences as elements of the space, and in that sense it has to be a "closed" space.

## 3.6   Bolzano-Weierstrass theorem

The Bolzano-Weierstrass theorem is a mathematical statement that, on first glance, may seem strange (and even appear wrong), and whose significance and importance might not be obvious. It is a core result in mathematical analysis; it is widely used in proofs in functional analysis and elsewhere, and we will require it below. In its simplest form, it can be stated as follows:

$$\text{Each bounded sequence of real numbers has a convergent subsequence} \tag{7}$$

Note, first, that the "sequence" here is, again, presumed to be infinitely long, even though this is typically not explicitly stated (and without this presumption the statement would indeed be wrong). Second, the bounded set of real numbers considered here has to include the supremum and infimum of the sequence – that is, it has to be a "closed" set of real numbers; otherwise, the subsequence could converge to a number that lies ("just") outside the set, which would make it technically not a convergent sequence. Third, the subsequence in the statement of the theorem is itself a sequence, and is also presumed to be an infinitely long sequence (which again is not typically explicitly stated).

We can consider a subsequence to be constructed from the original sequence by selecting some (or possibly all) of the elements while keeping their relative order. So, with a sequence $(1.1, 3.2, 2.5, 0.7, \ldots)$, $(3.2, 0.7 \ldots)$ is a possible subsequence, but $(0.7, 3.2 \ldots)$ is not a possible subsequence because the relative order of the elements has been changed[14].

---

[12] Cauchy sequences are also sometimes called *fundamental sequences*.

[13] Unfortunately, the term "complete" is used for more than one different purpose in this field. Here we are explicitly discussing a "complete space". Elsewhere, we will discuss a "complete set", which refers to a set of basis functions or vectors that can be used in linear combinations to make up any function or vector in a given space.

[14] For the notation "…", we may use this especially when specifying elements of a set or a sequence. When we write a set in the form $\{x_1, x_2, \ldots, x_n\}$ or a sequence in the form $(x_1, x_2, \ldots, x_n)$, we mean that there is some finite number, $n$, of elements in the set or sequence, and the "…" indicates we should include all the elements, continuing in the pattern given by the first few (here two) stated. So here we should be including the elements $x_3, x_4, x_5$ and so on (presuming here the $n > 5$), which, in the case of the sequence, should be in this order. When we write $\{x_1, x_2, \ldots\}$ or $(x_1, x_2, \ldots)$, we similarly mean that the set or sequence should continue with the next elements in the obvious pattern, but either we are presuming the set or sequence is infinite (which will be more common) or that it might be either finite (with a number of elements we are not specifying) or infinite.

It may help in understanding this theorem to look at some extreme cases. First, it is, of course, possible to construct an infinitely long sequence that does not itself converge. A simple example would be an "oscillating" sequence, such as $(1,0,1,0,\ldots)$. But this sequence does have two obvious convergent subsequences – explicitly $(1,1,1,\ldots)$ and $(0,0,0,\ldots)$. Each of these trivially converges, the first to 1, and the second to 0.

A point or value at which an infinite subsequence converges is called an *accumulation point*. An accumulation point is not necessarily a limit of the original sequence, which, after all, may not even have a limit to which it converges, but it is the limit of at least one subsequence. Obviously, in our oscillating sequence, there are two accumulation points, one being 1 (the limit of the subsequence $(1,1,1,\ldots)$) and the other being 0 (the limit of the subsequence $(0,0,0,\ldots)$) Our original oscillating sequence has many other convergent subsequences (in fact, in infinite number) because it is only necessary that the subsequence eventually converges; the "non-converging" part of it can go on as long as we like provide the subsequence eventually does converge. So, a subsequence $(1,0,1,0,0,0,\ldots)$, where all the remaining elements are 0, is also a convergent subsequence.

One other key point to note is that, in constructing a sequence in the first place, we can repeat the same element of the set as many times as we like (as we have done in constructing the oscillating sequences above). So, trivially, we can always construct a convergent sequence from any (non-empty) set of real numbers by just repeating the same number infinitely. This rather trivial kind of sequence is sometimes used in proofs.

We give one of the standard proofs of this theorem (7) below 11.2 "Proof (2) of the Bolzano-Weierstrass theorem".

## 3.7 Compact spaces

The formal definition of a *compact* space is as follows:

> A metric space is said to be *compact* if every sequence in the space has a convergent subsequence. (8)

Note that the convergence here will be in terms of the metric, as in Eq. (4) or Eq. (5), which is always defined in terms of the distance between two elements. Convergence here is *not* just convergence of the norm, which would be the "length" of a given element. For real numbers we might not notice the difference between a convergence in the metric and a convergence in the norm, but for other entities such as functions or vectors, we need to be considering convergence in the metric, and unless otherwise stated, that is what we mean. (Different vectors could all have the same norm, but be pointing in different "directions", in which case a sequence of them could converge with respect to the norm but not with respect to the metric.)

For a space to be compact, it is important that it is bounded; if the space is not bounded then we do not know if sequences in general have convergent subsequences; the Bolzano-Weierstrass theorem is not then applicable (it is easy to construct a sequence of numbers that just keeps on growing and hence has no convergent subsequence). But, if we have a complete space, so every convergent sequence converges within the space (including Cauchy convergence), then we can use the Bolzano-Weierstrass theorem to conclude that such a bounded space is compact. So, we expect bounded complete spaces to be compact.

A space is called *precompact* (or, in some texts "relatively compact") if its closure is compact. So a precompact space would be one in which some sequences might converge "just" outside the space, and the closure of that space, adding in the limits of such sequences to the space, would be compact.

# 4 Convergence in vector and function spaces

So far, we have explicitly discussed mostly just real numbers as the elements in a space, with suggestions that we are going to generalize beyond this. Now we are going to explicitly introduce the kinds of elements that are of most interest to us, which are vectors and functions. Generally, for us, these two concepts of vectors and functions are interchangeable. All of our vectors can be considered to be functions and all of our functions can be considered to be vectors. The general mathematics below will be the same for vectors and functions. (The only substantial difference will come in the precise way we choose to define what we will call the "inner product" below.) Because of this similarity, we can use the same notation for both.

## 4.1 Notation for scalars, vectors and functions

In mathematics texts, it is common to make no substantial distinction between the notation for a scalar or for a vector or function; we might use $x$ or $y$ for any of these, for example. From now on, we will try to be more specific, and, as much as possible, we use specific italic Roman letters, such as $a$, $b$, $c$, $d$, $f$, $g$, $r$, $s$, $t$, $x$, $y$, and $z$, for real or complex scalars, other italic Roman letters, such as $j$, $k$, $m$, $n$, $p$, and $q$, for positive integers, and italic Greek letters, such as $\alpha$, $\beta$, $\gamma$, $\mu$, $\eta$, $\theta$, $\phi$, and $\psi$, to indicate functions or vectors (though we may have to make occasional exceptions to this notation).

Later, we will be able to use Dirac's bra-ket notation more extensively for functions and vectors, such as the "ket" $|\phi\rangle$; for the moment, we can introduce this as an alternative notation for vectors, postponing for the moment the broader transition to using it also for functions.

By "vector" here, we usually mean something that we could also write out as a mathematical column vector of scalars, as in

$$\phi = \begin{bmatrix} f_1 \\ f_2 \\ f_3 \\ \vdots \end{bmatrix} \equiv |\phi\rangle \tag{9}$$

Note that we do not in general mean a geometrical vector – i.e., something with "components" along geometrical axes (though the mathematics can handle that as a specific case). If necessary to make the distinction, we will use the words "mathematical" on the one hand and "geometrical" or "physical" on the other hand, and by default we mean the "mathematical" version here, as in Eq. (9).

This mathematical vector could have a finite number of elements (and hence be of a finite "dimension"), though generally we want to set up the mathematics to handle vectors with possibly infinite numbers of elements (and hence of possibly infinite dimension); vectors with finite numbers of elements are then just special cases of vectors with generally infinite dimension.

In a sense, such a vector is already a function; it is a mapping from the natural numbers $1, 2, 3, \dots$ (the subscripts in the individual elements in the vector in Eq. (9)) to the corresponding values $f_1, f_2, f_3, \dots$. Once we define the inner product below, we will have a particularly definite way to clarify the relation between vectors and functions.

By a function here, we could mean a simple scalar function of one continuous real variable $x$, such as

$$\phi \equiv \phi(x) \tag{10}$$

We might also mean a function of multiple variables, as in a scalar function of the ordinary Cartesian $(x, y, z)$ coordinates, such as some function $\gamma(x, y, z)$, and we could even mean a physical vector function such as some electric field $\mathbf{E}(x, y, z) \equiv \mathbf{E}(\mathbf{r})$ (the use of the bold $\mathbf{E}$ for the function is an obvious exception to our lower case Greek letters for functions).

## 4.2 Vector (or function) space

In what follows, we will use the term "vector space" to cover spaces with either vectors or functions as their elements. We need to give such spaces some properties so we can work towards ideas of convergence. We can formally define a *vector (or function) space* (also known as a *linear space*) as a (non-empty) set contain vector (or function) elements such as $\alpha$, $\beta$, and $\gamma$, and having two algebraic operations: *vector addition* and *multiplication of vectors by scalars*. By scalars here we will mean complex numbers[15]. For the vectors and functions we are considering, such additions of vectors and of functions (which are just point-by-point or element-by-element additions) and multiplications of vectors or functions by scalars are both relatively obvious, so we defer these formal points to a footnote[16].

## 4.3  Matrices, vectors, Hermitian adjoints and Dirac notation

Though we will not need matrices until we have introduced linear operators, we do need to introduce some more notation for vectors, and it will save some space if we formally introduce all the basic notation and operations of matrices and vectors here, including more aspects of the Dirac notation.

First, we should explicitly clarify what we mean by a matrix (even though this is very standard). A *matrix* is a 2-dimensional array of numbers (which we will allow to be complex). When using a letter to refer to a matrix, we will use an upper-case sans-serif font, as in A. (We will also below use this same convention when referring generally to linear operators, which we will define later.) By convention, the elements of a matrix are written with two subscripts, in the order of the "row number" first and the "column" number second, where the numbering starts in the top left corner, and uses the natural numbers (so, starting at 1, not 0). We will typically use the lower case italic Roman letter corresponding to the letter name of the matrix, so $a_{jk}$ for the element in the $j$th row and $k$th column of the matrix A. Hence the matrix A can be written as

$$A = \begin{bmatrix} a_{11} & a_{12} & \dots \\ a_{21} & a_{22} & \cdots \\ \vdots & \vdots & \ddots \end{bmatrix} \tag{11}$$

If the matrix has finite dimensions, then we can characterize it by the number of rows, $N$, and the number of columns, $M$, saying such a matrix is an $N \times M$ matrix.

With two vectors written as

---

[15] Technically, a vector space can be defined with respect to any specific mathematical "field", but we will exclusively be considering the field of complex numbers, of which real numbers are also just a special case.

[16] Vector addition is an operation that is both commutative, i.e., $\alpha + \beta = \beta + \alpha$, and associative, i.e., $\alpha + (\beta + \gamma) = (\alpha + \beta) + \gamma$. To deal fully with convergence, we require that the space needs to include a "zero" vector. We could write such a vector as, for example, *zero* to make a distinction in our notation; however, generally both mathematicians and physical scientists are very loose here, and general just write "0" instead for this vector, on the presumption that no-one will be confused into thinking this is a number rather than a vector. With that dubious convention, we formally have $\alpha + zero \equiv \alpha + 0 = \alpha$ for all vectors $\alpha$. Also, we require that for every vector $\alpha$ in the space, there is a vector $-\alpha$ such that $\alpha + (-\alpha) = 0$ (where the "0" on the right is the zero vector, not the number zero). For all vectors $\alpha$ and $\beta$ and all scalars $a$ and $b$, multiplication by scalars, usually written in the form $a\alpha$, is such that we have $a(b\alpha) = (ab)\alpha$ and $1\alpha = \alpha$ (where the "1" here is the real number 1 – the multiplicative identity element in the field of complex numbers). The usefulness of this is in complicated multiplicative vector expressions (and we will define what we mean by vector multiplications later), where we note that we can move scalars essentially at will through any such products. Note too that the multiplication by scalars is effectively commutative – we can write $a\alpha = \alpha a$ (even though we typically use the first of these notations). We also have the distributive laws $a(\alpha + \beta) = a\alpha + a\beta$, $\alpha(x + y) = \alpha x + \alpha y$ and $(\alpha + \beta)x = \alpha x + \beta x$. For the case when we are dealing with functions such as an electric field $\mathbf{E}(\mathbf{r})$ that is itself a (geometrical) vector-valued function, the addition of such functions should, of course, be (geometrical) vector addition, but such (geometrical) vector arithmetic operations obey the same formal rules as scalar functions with regard to associativity in addition and distributivity when multiplying by scalars.

$$\beta = \begin{bmatrix} b_1 \\ b_2 \\ \vdots \end{bmatrix} \equiv |\beta\rangle \quad \text{and} \quad \gamma = \begin{bmatrix} g_1 \\ g_2 \\ \vdots \end{bmatrix} \equiv |\gamma\rangle \tag{12}$$

we could have that the vector $|\gamma\rangle$ is the result of the multiplication of $|\beta\rangle$ by the matrix A, which we could write in any of the four equivalent ways:

$$\begin{bmatrix} g_1 \\ g_2 \\ \vdots \end{bmatrix} = \begin{bmatrix} a_{11} & a_{12} & \cdots \\ a_{21} & a_{22} & \cdots \\ \vdots & \vdots & \ddots \end{bmatrix} \begin{bmatrix} b_1 \\ b_2 \\ \vdots \end{bmatrix} \tag{13}$$

$$g_j = \sum_k a_{jk} b_k \tag{14}$$

$$\gamma = \mathsf{A}\beta \tag{15}$$

$$|\gamma\rangle = \mathsf{A}|\beta\rangle \tag{16}$$

where the summation form (14) is the most explicit about the actual details of the process of matrix-vector multiplication. For matrix-vector situations, we will typically prefer the bra-ket notation (16) over the more general notation (15).

We could, of course, have written a vector just as a special case of a matrix – a matrix with just one column – but in our use we have quite different physical meanings for vectors and matrices. Typically, a mathematical vector (not a geometric vector) will refer to some physical field, such as an acoustic or electromagnetic field, and a matrix will refer to a mathematical operation we perform on such a physical field or to some physical process such as one that generates waves from sources. As a result, it is useful for us to make an explicit distinction between vectors and matrices in our notation.

Because we want to work with complex vectors and functions (i.e., ones with complex-numbered values), we need a version of complex conjugation that is particularly convenient for the algebra of working with entire vectors and matrices, and this concept is the Hermitian adjoint[17]. For a vector (and also for a matrix), the *Hermitian adjoint* is formed by "reflecting" the vector or matrix about a "-45°" line (as in the "leading" diagonal of the matrix, with elements $a_{11}, a_{22}, \ldots$) – the operation known as "taking the transpose" – and then taking the complex conjugate of all the elements. This operation is usually notated with a superscript "dagger", "$\dagger$". So, for a matrix we have

$$\mathsf{A}^\dagger \equiv \begin{bmatrix} a_{11} & a_{12} & \cdots \\ a_{21} & a_{22} & \cdots \\ \vdots & \vdots & \ddots \end{bmatrix}^\dagger \equiv \begin{bmatrix} a_{11}^* & a_{21}^* & \cdots \\ a_{12}^* & a_{22}^* & \cdots \\ \vdots & \vdots & \ddots \end{bmatrix} \tag{17}$$

where the superscript "$*$" denotes the complex conjugate of the number.

The Hermitian adjoint of a (column) vector is a row vector with complex conjugated elements. Because we use this operation often with vectors in the algebra that follows, it also has its own notation, which is the "bra" part, $\langle\beta|$, of Dirac's bra-ket notation. Explicitly,

$$|\beta\rangle^\dagger \equiv \begin{bmatrix} b_1 \\ b_2 \\ \vdots \end{bmatrix}^\dagger \equiv \begin{bmatrix} b_1^* & b_2^* & \cdots \end{bmatrix} \equiv \langle\beta| \tag{18}$$

Note that in general the Hermitian adjoint operation performed twice gets us back to where we started, i.e., $\left(\mathsf{A}^\dagger\right)^\dagger \equiv \mathsf{A}$, and $\left(|\beta\rangle^\dagger\right)^\dagger \equiv \left(\langle\beta|\right)^\dagger \equiv |\beta\rangle$.

---

[17] The Hermitian adjoint is also known as the *Hermitian conjugate, conjugate transpose* or sometimes just the *adjoint*.

## 4.4    Inner products

For the vectors or functions of interest to us for our physical problems, we will be able to define what we can call an *inner product*. An inner product has a specific set of properties we will define below, but one of its most important characteristics is that it is a kind of product between two vectors or functions $\beta$ and $\gamma$ that results in a scalar number $c$ as the result. One mathematical notation for this is to write

$$c = (\beta, \gamma) \tag{19}$$

Before giving the formal properties of an inner product, we can give some simple examples.

### 4.4.1    Cartesian inner product

For (mathematical) vectors, a particularly simple inner product (but not the only one we could define) is what can be called the *Cartesian inner product*, which is just the vector product of two vectors constructed by multiplying one vector by the Hermitian adjoint of the other so as to give a scalar (number) result, as in

$$(\beta, \gamma)_{Cartesian} \equiv \begin{bmatrix} b_1^* & b_2^* & \cdots \end{bmatrix} \begin{bmatrix} g_1 \\ g_2 \\ \vdots \end{bmatrix} \equiv \sum_j b_j^* g_j \equiv \langle \beta \| \gamma \rangle \equiv \langle \beta | \gamma \rangle \tag{20}$$

On the right hand side, we show another standard part of the Dirac notation, which is to condense this into the compact notation $\langle \beta | \gamma \rangle$ that is a shorthand for $\langle \beta \| \gamma \rangle$. In our discussion here, we are going to use the Dirac notation for all such Cartesian inner products, though we will keep the general notation $(\beta, \gamma)$ also for other inner products or for the general case.

Another simple example of an inner product would be an integral form for functions of one variable, for example,

$$\langle \alpha | \beta \rangle \equiv \int \alpha^*(x) \beta(x) dx \tag{21}$$

This particular integral inner product is essentially also a Cartesian inner product in the limit of a sum tending to an integral, and so we have used the Dirac notation.

An integral like Eq. (21) is also a good example of what can be called a "functional". A mathematical definition of a *functional* is a mapping from a vector or function space to a space of scalars. In other words, a functional turns a vector or function into a number (just as an operator turns a vector or function into (usually another) vector or function). We could view the operation in Eq. (21) as a functional acting on the function $\beta(x)$ to generate the number $\langle \alpha | \beta \rangle$ on the left. In this article, the only functional we need is the inner product, so from this point on, we will not discuss functionals[18] in general any further[19].

An even simpler example of a Cartesian inner product is the usual (geometrical) vector dot product for geometrical vectors **a** and **b**, so we could write

---

[18] That we are avoiding any general treatment of functionals might well be considered by some to be almost an indictable offense in an introduction to functional analysis. Functionals generally have other important uses and were very important in the development of this field, which largely grew out of the need to solve integral equations. However, the most important result from the theory of functionals for our purposes is the inner product, and it is indeed very powerful. Our omission of any more general discussion of functionals is also one of the ways we can keep this introduction short and to the point.

[19] One other important example is a Green's function equation of the form $g(x_2) = \int G(x_2; x_1) f(x_1) dx_1$; for example, we might have a Green's function $G(x_2; x_1) \propto \exp(ik|x_2 - x_1|) / |x_2 - x_1|$ for a scalar wave equation, with $g(x_2)$ being the wave generated at $x_2$ from the source function $f(x_1)$. With $x_2$ viewed as a parameter, then this integral would be a functional, with $g(x_2)$ just being a number. However, we prefer to think of this as an integral operator $\int G(x_2; x_1) dx_1$ acting on the function $f(x_1)$ to generate the new function $g(x_2)$, so we do not use the "functional" way of looking at such an equation.

$$\langle \mathbf{a} | \mathbf{b} \rangle = \mathbf{a} \cdot \mathbf{b} \tag{22}$$

One main difference in inner products compared to the geometrical vector dot product is that the general inner product is designed for complex numbered elements, and that means that the order of the inner product generally matters (see (IP3) below).

### 4.4.2   Formal properties of inner products

For subsequent uses of inner products in the mathematics that follows, we only need the inner product to have the following properties, and this gives the most general definition[20].

For all vectors $\alpha$, $\beta$ and $\gamma$ in a vector space, and all (complex) scalars $a$, we define an inner product $(\alpha, \beta)$, which is a (complex) scalar, through the following properties:

(IP1)      $(\gamma, \alpha + \beta) \equiv (\gamma, \alpha) + (\gamma, \beta)$

(IP2)      $(\gamma, a\alpha) = a(\gamma, \alpha)$ (where $a\alpha$ is the vector or function in which all the values in the vector or function $\alpha$ are multiplied by the (complex) scalar $a$)

(IP3)      $(\beta, \alpha) = (\alpha, \beta)^*$

(IP4)      $(\alpha, \alpha) \geq 0$, with $(\alpha, \alpha) = 0$ if and only if $\alpha = 0$ (the zero vector)

$$\tag{23}$$

For the specific case of the Cartesian inner product, these criteria can be rewritten as

(IP1)      (Cartesian)      $\langle \gamma | \alpha + \beta \rangle \equiv \langle \gamma | \big( |\alpha\rangle + |\beta\rangle \big) \rangle = \langle \gamma | \alpha \rangle + \langle \gamma | \beta \rangle$

(IP2)      (Cartesian)      $\langle \gamma | a\alpha \rangle = a \langle \gamma | a \rangle$ (where $|a\alpha\rangle$ is the vector or function in which all the values in the vector or function $|\alpha\rangle$ are multiplied by the (complex) scalar $a$)

(IP3)      (Cartesian)      $\langle \beta | \alpha \rangle = \langle \alpha | \beta \rangle^*$

(IP4)      (Cartesian)      $\langle \alpha | \alpha \rangle \geq 0$, with $\langle \alpha | \alpha \rangle = 0$ if and only if $|\alpha\rangle = 0$ (the zero vector)

$$\tag{24}$$

We can easily check that all of our examples above of inner products satisfy all these four criteria, either in the form (23) or the Cartesian form (24).

Note that, in addition to (IP2) above, we can write as a consequence of these properties that

$$\begin{aligned}
\left(a\gamma,\alpha\right) &= \left(\alpha,a\gamma\right)^* \text{ by (IP3)} \\
&= \left[a\left(\alpha,\gamma\right)\right]^* \text{ by (IP2)} \\
&= a^*\left(\alpha,\gamma\right)^* \\
&= a^*\left(\gamma,\alpha\right) \text{ by (IP3)}
\end{aligned} \tag{25}$$

---

[20] Note, incidentally, that the order in which we are writing the inner product here is the opposite order from that used in most mathematics texts. This difference shows up specifically in (IP2); most mathematics texts would write $(a\alpha, \gamma) = a(\alpha, \gamma)$ instead. The convention in these mathematics texts is unfortunate because it is the opposite way round from the order we find in matrix-vector multiplication (and in Dirac notation). The matrix-vector notation allows a simple associative law without having to change the written order of the elements, whereas this conventional mathematics notation does not. The Dirac notation follows the matrix-vector ordering (and indeed, Dirac notation is generally a good notation for complex matrix-vector algebra). At least one modern text [3] recognizes the problems of this historical choice in mathematics texts, and uses the notation $(\alpha, \gamma)$ the other way round, as we do here.

The combination of properties $\left(\gamma,a\alpha\right)=a\left(\gamma,\alpha\right)$ and $\left(a\gamma,\alpha\right)=a^*\left(\gamma,\alpha\right)$ means the inner product is what is sometimes called *sesquilinear*. "Sesqui" is the Latin word for "one and a half", loosely indicating that the inner product is only "half" linear when the factor is in front of the left vector because we then require the complex conjugate of the multiplying factor.

### 4.4.3   Weighted inner products

We might ask what additional possibilities are created by the more general definition of the inner product in (23) rather than just the Cartesian one in (24) above. There is one particularly important extension that this general definition allows, which is to what we could call a *weighted inner product*. An example of a weighted inner product might be

$$\left(\alpha,\beta\right)=\int\gamma(x)\alpha^*(x)\beta(x)dx \tag{26}$$

where $\gamma(x)$ is a real, positive, and non-zero function of $x$. This is an inner product in the sense of satisfying all of the general criteria in (24). Such weighted inner products can occur in physical problems.

An example would be in an inner product that can give the electrostatic energy corresponding to a field $\mathbf{E}(\mathbf{r})$ in a dielectric with a scalar, positive, non-zero dielectric constant $\varepsilon(\mathbf{r})$. Then the dielectric constant $\varepsilon(\mathbf{r})$ (or $\varepsilon(\mathbf{r})/2$) would be the weighting function, and we could define the inner product (here for the specific case of an inner product of a field with itself) as

$$W=\left(\mathbf{E},\mathbf{E}\right)=\iiint\frac{1}{2}\varepsilon(\mathbf{r})\mathbf{E}^*(\mathbf{r})\cdot\mathbf{E}(\mathbf{r})d^3\mathbf{r} \tag{27}$$

Even with the presence of the weighting function and with the (physical) vector dot product as the multiplication, this satisfies all the criteria in (23) above, and so is a valid inner product. It is also an example of an *"energy" inner product*, where the inner product of a vector or function with itself gives the energy $W$ in the system.

Below, in section 8, we show a further broad category of entities that also are inner products, though we defer this discussion for the moment until we have introduced Hermitian operators.

## 4.5   Norms and metrics from inner products

We note that the inner product of a vector or function with itself is guaranteed to be zero or a positive real number[21]; this is the criterion (IP4) above in (23). As a result, the inner product can be used to define a norm for a vector or function $\alpha$, which we can choose as

$$\|\alpha\|=\sqrt{\left(\alpha,\alpha\right)} \tag{28}$$

which we can think of as being the "length" or the magnitude of the "amplitude" of the vector or function. (This norm satisfies all the criteria for a norm as in (1).)

As is generally true for norms, they can be used to define a metric. So, for some vector space $P$ in which these vectors are elements, we can therefore define the metric

---

[21] Since an inner product must satisfy $\left(\beta,\alpha\right)=\left(\alpha,\beta\right)^*$, which is criterion (IP3) in (23), $\left(\alpha,\alpha\right)$ is guaranteed to be a real number since it must equal its own complex conjugate.

$$d_P(\alpha, \beta) \equiv \|\alpha - \beta\| = \sqrt{(\alpha - \beta, \alpha - \beta)} \tag{29}$$

A (vector or function) space with an inner product defined on it can be called an *inner product space*. So, all inner product spaces are normed spaces (and also, of course, metric spaces).

## 4.6    Convergence

Our previous mathematical arguments on convergence were nominally written for convergences of sequences of real numbers. Because we wrote them in terms of a metric, we can, however, now extend those same arguments and definitions, without change, to vector or function spaces; we just substitute our new metric as in Eq. (29), and consider the elements in the space to be vectors like $\alpha$ and $\beta$ instead of real numbers $x$ and $y$. We can therefore talk about a sequence of vectors, which we could write as $(\alpha_n)$, and consider convergent sequences, including Cauchy sequences, of vectors, and we can have complete vector spaces in the same sense as complete spaces of real numbers.

## 4.7    Inner products and orthogonality

The other very important property of an inner product is that we use it to define the concept of *orthogonality*. Specifically,

a non-zero element $\alpha$ of an inner product space is said to be *orthogonal* to a non-zero element $\beta$ of the same space if and only if

$$(\alpha, \beta) = 0 \tag{30}$$

This is a generalization of the idea of the geometrical vector "dot" product, which is similarly used to define orthogonality in geometrical space. Note here that we are extending that idea to allow for complex vector "components" and for arbitrary, even possibly infinite, numbers of dimensions.

## 4.8    Energy inner products and orthogonality

Such orthogonality as defined using the inner product, in addition to being very useful as a mathematical device, also has direct physical meaning and implications. For example, suppose we are able to write some electrostatic field $\mathbf{E}$ as the sum of two (non-zero) fields $\mathbf{E}_1$ and $\mathbf{E}_2$ that are mathematically orthogonal as defined by their energy inner product, i.e., $(\mathbf{E}_1, \mathbf{E}_2) = 0$, as in Eq. (27). Then, the energy associated with this field would be

$$W = (\mathbf{E}, \mathbf{E}) = (\mathbf{E}_1 + \mathbf{E}_2, \mathbf{E}_1 + \mathbf{E}_2) = (\mathbf{E}_1, \mathbf{E}_1) + (\mathbf{E}_2, \mathbf{E}_2) + (\mathbf{E}_1, \mathbf{E}_2) + (\mathbf{E}_2, \mathbf{E}_1)$$
$$= (\mathbf{E}_1, \mathbf{E}_1) + (\mathbf{E}_2, \mathbf{E}_2) \tag{31}$$

So, if the two fields are orthogonal in their energy inner product, then the total energy is just the sum of the energies of each field considered separately. Such energy inner products are therefore particularly important for physical fields[22].

---

[22] Furthermore, if we want to quantize the fields, it is very desirable to start with classical fields that are orthogonal in such an energy inner product; then we can separate the Hamiltonian into a sum of Hamiltonians, one for each field that is orthogonal to all the others, where orthogonality is defined using this inner product, and then quantize those Hamiltonians separately.

# 5  Hilbert spaces and their properties

The single most important kind of vector space for our purposes is a Hilbert space, and we will be working with these from this point on. With the mathematics we have constructed up to this point, the definition of a *Hilbert space* is now straightforward[23].

$$\text{A Hilbert space is a complete inner product space.} \tag{32}$$

Generally,

> we will notate spaces (which might or might not be Hilbert spaces) using italic upper case letters such as *D, F, G,* and *R,* with subscripts to distinguish spaces if necessary,

> and if we use *H* (with or without a subscript) for a space, then it is certainly a Hilbert space.

We have already seen that inner product spaces have physical meaning and usefulness. The notion of completeness is generally something we can expect our spaces of physical vectors to have. Essentially it means we are not "missing out" specific functions from the space, and we are careful to include functions that might be at the "edges" in a mathematical sense of our space – i.e., that are the limits of sequences of functions that are otherwise within the space.

One of the reasons why Hilbert spaces are so useful is that we can construct orthogonal sets of functions, and use those to represent other functions. Indeed, once we construct appropriate sets, called basis sets, we can represent any function in the space using them.

## 5.1  Orthogonal sets and basis sets

An *orthogonal set* in a Hilbert space is a subset of the space whose (non-zero) elements are pairwise orthogonal (i.e., every member is orthogonal to every other member). So, for any two (non-zero) members $\alpha$ and $\gamma$ of this orthogonal set, $(\alpha, \gamma) = 0$ unless $\alpha = \gamma$. Even more convenient is an *orthonormal set*, which is an orthogonal set in which every element has norm 1, i.e., $(\alpha, \alpha) = 1$; functions with such a norm can be called *normalized*. We presume that we can index the members with an integer or natural number index[24], *j* or *k*, for example, in which case we can write for an orthonormal set

---

[23] In making this step to Hilbert spaces, we have "jumped over" Banach spaces. Banach spaces are complete metric (vector) spaces, but do not necessarily have an inner product. These are typically discussed at length in mathematics texts, but we have no real use for them here, so we omit them. Some of the definitions theorems and proofs we use in Hilbert spaces can be executed in the more general Banach spaces but anything that is true for inner product spaces in general or for Banach spaces is also true Hilbert spaces, because Hilbert spaces are just special cases of Banach spaces (i.e., with the explicit addition of the inner product). Similarly, a few definitions and results can be constructed in inner product spaces that are not necessarily complete, but Hilbert spaces are again just special cases of these. So here we just give all results for Hilbert spaces without calling out those that also apply to the simpler inner product or Banach spaces.

[24] Technically, this assumption that the orthogonal or orthonormal sets of interest to use can be indexed with integers or natural numbers, extending possibly to an infinite number, is an assumption that the set is *countable*. (A countable set is simply one whose elements can be put in one-to-one correspondence with members of the set of integers or natural numbers; countable sets include the set of all rational numbers, but the set of real numbers, for example, is not itself countable.) From this point on, we are technically presuming our Hilbert spaces are countable in this sense. We could argue that we can justify such an assumption *a posteriori* because our resulting mathematics works in modeling physical problems, which is in the end the only justification for using any mathematical models for the physical world. For physical problems, such an assumption of countability is common and implicit in constructing basis sets. For example, in working with plane wave functions in all directions (a set that is not countable if the components of the direction vectors are taken to be arbitrary real numbers), it is common to imagine a box of finite but large size, and count the functions that fit within the box, with "hard wall" or periodic boundary conditions at the edges of the box. This is an *ad hoc* construction of a countable set of plane wave functions, and it remains countable as the size of the box is increased towards infinity in the various directions.

$$\left(\alpha_j, \alpha_k\right) = \delta_{jk} \tag{33}$$

where the *Kronecker delta* $\delta_{jk}$ has the properties $\delta_{jk} = 1$ for $j = k$, but $\delta_{jk} = 0$ otherwise.

An important property of such sets is the notion of linear independence. To discuss this, we first formally need to introduce the idea of linear dependence and some related terms. A *linear combination* of vectors $\beta_1, \ldots, \beta_m$ of a vector space is an expression of the form $d_1\beta_1 + d_2\beta_2 + \ldots + d_m\beta_m$ for some set of scalars $\{d_1, d_2, \ldots, d_m\}$.

The idea of whether a set of vectors is linearly independent has to do with whether one (or more) members of the set can be expressed as linear combinations of the others; if this is possible, then the set is linearly dependent; if not, the set is linearly independent. Formally, this is decided using the equation $d_1\beta_1 + d_2\beta_2 + \ldots + d_m\beta_m = 0$. (We presume that there is at least one element here, i.e., $m \geq 1$.) If the only set of scalars $\{d_1, d_2, \ldots, d_m\}$ for which this holds is when they are all zero, then the set of vectors $\{\beta_1, \ldots, \beta_m\}$ is *linearly independent*. Otherwise, there is always a way of expressing some vector in the set in terms of a linear combination of the others, and the set of vectors is *linearly dependent*[25]. Note that an orthogonal or orthonormal set is linearly independent[26].

We can choose to have a set of vectors defined as those vectors $\gamma$ that can be represented in an orthonormal set $\{\alpha_1, \alpha_2, \ldots\}$ by the sum

$$\gamma = a_1\alpha_1 + a_2\alpha_2 + \cdots \equiv \sum_j a_j\alpha_j \tag{34}$$

We can also call such an expression the *expansion* of $\gamma$ in the basis $\alpha_j$ (i.e., in the set $\{\alpha_1, \alpha_2, \ldots\}$), and the numbers $a_j$ are called the *expansion coefficients*.

We can give the corresponding space of all such vectors that can be written using such an expansion the name[27] $G_\alpha$. This is then the set of all vectors that can be represented using this "basis" set $\{\alpha_1, \alpha_2, \ldots\}$. A set of orthogonal vectors (and, preferably, orthonormal vectors for convenience) that can be used to represent any vector in a space can be called an (orthogonal or orthonormal) *basis* for that space. Indeed, because we have deliberately constructed this set using only linear combinations of this set of orthonormal vectors $\{\alpha_1, \alpha_2, \ldots\}$, this set is automatically a basis for the space $G_\alpha$. The number of orthogonal or orthonormal functions in the basis – i.e., the number of functions in the set $\{\alpha_1, \alpha_2, \ldots\}$ – is called the *dimensionality* of the basis set and of the corresponding space. Depending on the space, this dimensionality could be finite or it could be infinite[28]. The basis set that can be used to construct any function in a given space is said to *span* the space.

---

[25] E.g., for non-zero $d_m$, we could write $-\left(1/d_m\right)\left(d_1\beta_1 + d_2\beta_2 + \ldots + d_{m-1}\beta_{m-1}\right) = \beta_m$. $\beta_m$ is then being expressed as a linear combination of the other vectors.

[26] To prove the linear independence of an orthogonal set formally, consider the orthogonal set of vectors $\{\alpha_1, \ldots, \alpha_n\}$ and consider the equation $d_1\alpha_1 + \cdots + d_n\alpha_n = 0$ (the zero vector) with the $d_j$ being complex numbers. Taking the inner product with any one of the elements $\alpha_j$ leads to $\left(\alpha_j, \left(d_1\alpha_1 + \cdots + d_n\alpha_n\right)\right) = d_j\left(\alpha_j, \alpha_j\right) = 0$ (the number zero). Since the elements $\alpha_j$ are by definition non-zero, this implies that every $d_j$ is zero, which means that the set is linearly independent (no vector in this set can be made up from a linear superposition of other vectors in the set).

[27] We have not yet proved that this space formed in this way by such vectors is a Hilbert space, though we can always find such sets of orthogonal functions for a Hilbert space, as is proved later.

[28] Indeed, much of our reason for setting up this formal mathematics is because we need to deal with spaces of possibly infinite dimensionality. If we were only dealing with finite dimensionality, the mathematics can be expressed more simply, but we need the infinite-dimensional results. The results for finite dimensionalities are then just special cases.

By definition, a basis, because it can represent any function in a given space, is also said to be a *complete set* of functions for the space. (Note, incidentally that this is a different use of the word "complete" from the idea of a complete space as defined above; this potential confusion is unfortunate, but is unavoidable in practice because of common usage[29].)

The coefficient $a_j$ in the expansion Eq. (34) is easily extracted by forming the inner product with $\alpha_j$, i.e.,

$$\left(\alpha_j, \gamma\right) = a_j \tag{35}$$

Indeed, we can take this to be the defining equation for the expansion coefficients. Note this evaluation of the coefficients uses whatever we have defined for the inner product of the space; the inner product in the Hilbert space of interest need not be a Cartesian inner product.

We can now view this set of numbers $\left\{a_j\right\}$ as being the *representation* of the function $\gamma$ in the basis $\left\{\alpha_1, \alpha_2, \ldots\right\}$, and we can choose to write them as a column vector

$$\gamma = \begin{bmatrix} a_1 \\ a_2 \\ \vdots \end{bmatrix} \equiv \left| \gamma \right\rangle \tag{36}$$

Now, a key conceptual point is that, since the function $\gamma$ will typically be an actual physical function – such as an electromagnetic field, for example – it is the same function no matter how we represent it. There are many basis sets (in fact, usually an infinite number of possibilities[30]) that can be used to represent a given function in a space, but no matter which representation we use, and therefore which specific column of numbers we are writing in an expression like Eq. (36), the function is the same function. (Indeed, this could be regarded as one justification why Dirac notation does not include any explicit specification of the basis – at one level, it makes no difference what the basis is.)

We should note explicitly that the specific form of the inner product is something that goes along with the Hilbert space and is part of the definition of the space. Indeed, it will be useful to give a name to the inner product used in the definition of a given Hilbert space; we will call it the *underlying inner product* of the Hilbert space[31]. Whatever basis we decide to use, its orthogonality should be set using this underlying inner product, and the expansion coefficients should be calculated using this underlying inner product.

## 5.2    Basis sets on Hilbert spaces

A key attribute of a Hilbert space is that there is always a basis for a Hilbert space. That is, there is always some complete set of orthogonal (or orthonormal) functions $\left\{\alpha_1, \alpha_2, \ldots\right\}$ that forms a basis for any given Hilbert space. Importantly for our purposes, this also applies to infinite-dimensional Hilbert spaces. We give the proof below in 11.3 "Proof (3) of the existence of a basis for a Hilbert space".

## 5.3    An "algebraic shift"

The existence of a basis set for any given Hilbert space is very useful when we describe linear operators, which is the subject of the next section, and it allows an important "algebraic shift", as we will now explain.

---

[29] Mathematics texts sometimes use the terminology *total set* instead of "complete set", but this is not common in physical science and engineering.

[30] Any linear combination of the original orthogonal or orthonormal basis sets that results in new orthogonal vectors can be used as a basis, and there is an infinite number of such possibilities.

[31] This terminology of "underlying inner product" is one that I am introducing here. It is not a standard term as far as I am aware.

Because we can always construct a basis set in a Hilbert space, then we can always construct a mathematical column vector as in Eq. (36) to represent an arbitrary function in the Hilbert space. This means that, once we have constructed the expansion coefficients using the underlying inner product in the space, as in Eq. (35), the subsequent inner products of functions represented with vectors of such expansion coefficients can always be considered to be in the simple "row-vector times column-vector" Cartesian form as in Eq. (20). It is this option to change subsequently to such Cartesian inner products that we are calling our *algebraic shift*.

To see why this works we can formally consider an inner product of two functions $\eta$ and $\mu$ in a given Hilbert space. To start, we expand each function in an orthonormal basis set $\{\alpha_1, \alpha_2, \ldots\}$ for the space, obtaining

$$\eta = \sum_k r_k \alpha_k \text{ and } \mu = \sum_k t_k \alpha_k \tag{37}$$

where

$$r_k = (\alpha_k, \eta) \text{ and } t_k = (\alpha_k, \mu) \tag{38}$$

are inner products formed using the underlying inner product in the space, which might be, for example, a weighted inner product such as an energy inner product.

Now the inner product of $\eta$ and $\mu$ in this Hilbert space can be written

$$(\mu, \eta) = \sum_{p,q} t_p^* r_q (\alpha_p, \alpha_q) = \sum_{p,q} t_p^* r_q \delta_{pq} = \sum_{p,q} t_p^* r_q \equiv \begin{bmatrix} t_1^*, t_2^*, \ldots \end{bmatrix} \begin{bmatrix} r_1 \\ r_2 \\ \vdots \end{bmatrix} \tag{39}$$

So, once we have made the "algebraic shift" of regarding the vectors as being vectors of expansion coefficients that have been constructed using the underlying inner product in the space, then the subsequent mathematics of the inner products is simply the Cartesian inner product as in Eq. (20).

So, because there is always an orthonormal basis for any Hilbert space, now we can always write any vector or function $\eta$ in a Hilbert space as the "ket" $|\eta\rangle$. We can consider this ket to be the column vector of numbers

$$|\eta\rangle = \begin{bmatrix} (\alpha_1, \eta) \\ (\alpha_2, \eta) \\ \vdots \end{bmatrix} \tag{40}$$

With the understanding, as in Eq. (18), that we can similarly write the Hermitian adjoint of any such vector as

$$\langle \eta | \equiv |\eta\rangle^\dagger = \begin{bmatrix} (\alpha_1, \eta) \\ (\alpha_2, \eta) \\ \vdots \end{bmatrix}^\dagger = \begin{bmatrix} (\alpha_1, \eta)^* & (\alpha_2, \eta)^* & \cdots \end{bmatrix} \tag{41}$$

then we can write the inner product of any two vectors $\mu$ and $\eta$ in a given Hilbert space as

$$(\mu, \eta) \equiv \begin{bmatrix} (\alpha_1, \mu)^* & (\alpha_2, \mu)^* & \cdots \end{bmatrix} \begin{bmatrix} (\alpha_1, \eta) \\ (\alpha_2, \eta) \\ \vdots \end{bmatrix} \equiv \langle \mu | |\eta\rangle \equiv \langle \mu | \eta \rangle \tag{42}$$

So, the inner product of any two vectors or functions in the space – an inner product that must be formed using the underlying inner product of the space – can be rewritten as a Cartesian inner product of the two vectors consisting of the expansion coefficients on a basis, where those expansion coefficients are formed using the underlying inner product.

Because we have now found a way of writing any inner product as a Cartesian inner product (sitting "above" the underlying inner product in the expansion coefficients), algebraically we can now "break up" the inner product into the simple "Cartesian" product of two vectors as in

$$\langle \mu | \eta \rangle \equiv \langle \mu \| \eta \rangle \tag{43}$$

even when the underlying inner product would not necessarily allow us to do this. This "algebraic shift" then allows us to use the full algebraic power of vector-matrix multiplication, including associative laws that break up the inner product as in Eq. (43). We will return to this once we have similarly considered representing operators in a related way.

This algebraic shift also gives us a specific way of seeing vectors as "being" functions: we can write out any function in the Hilbert space as such a mathematical column vector by performing the expansion using the underlying inner product.

# 6  Linear operators

An *operator* is something that turns one function into another, or, equivalently, generates a second function starting from a first function. Generally, an operator maps from functions in its *domain*, a space $D$, to functions in its *range*, a space $R$. Here, we will consider both the domain and the range to be Hilbert spaces. They may be the same space or they may be different spaces[32].

In our case, we are specifically interested in linear operators[33]. With a linear operator A, we write the action of the operator on any vector or function $\alpha$ in its domain $D$ to generate a vector or function $\gamma$ in its range $R$ as

$$\gamma = A\alpha \tag{44}$$

The linear superposition requirement is consistent with the usual definition of scalars and linear operators:

For any two vectors or functions $\alpha$ and $\beta$ in its domain $D$, and any scalar $c$ (which here we allow to be any complex number), an operator is a *linear operator* if and only if it satisfies the two properties:

(O1)     $A(\alpha + \beta) = A\alpha + A\beta$

(O2)     $A(c\alpha) = cA\alpha$

$$\tag{45}$$

In words, the first property, O1, says that we can calculate the effect of the operator on the sum of two vectors or functions by calculating its effect on the two functions separately and adding the result. The second property, O2, says that the effect of the operator on $c$ times a vector or function is the same as $c$ times the effect of the operator on the function.

---

[32] An example physical problem where the domain and range are quite different spaces is where we start with source functions in one volume that lead to waves in another volume. Not only would the generated functions be in a different space – actually, even a different physical volume – than the source functions; they could be built from entirely different physical quantities. The source functions might be current densities, and the resulting waves might be electric fields. We might therefore have quite different kinds of inner products in these two spaces. Situations like these could, however, be handled with operators mapping between the spaces. In our mathematics, we can also formally keep track of just what inner product is being used in each space; the underlying mathematics supports this even if it is not commonly explicit to have different inner products in different spaces.

[33] For example, we are presuming here that any physical wave systems we are considering are linear, with linear superposition applying to waves and sources.

## 6.1    Bounded operators and operator norms

To consider the convergence of the effects of operators, and even to consider the idea of convergences of sets of operators themselves, we need to define a *operator norm* for any operator A. There is more than one way of doing this, though like any other norm, any such operator norm must satisfy the properties given in (1) above for a norm. The first approach, which is quite general in that it can apply to a large range of types of operators, is what we can call a supremum (operator) norm. This norm is built on the norms for vectors or functions in both the domain and the range.

In this case, we consider the set of all possible (non-zero) vectors or functions $\alpha$ in the domain $D$. Any particular vector will have some norm $\|\alpha\|$ in the domain space (which we will take to be the one based on the underlying inner product in that space as in Eq. (28)). In the range space $R$, the vector A$\alpha$ will have a (vector) norm $\|A\alpha\|$ (based on the underlying inner product in that space), and we can consider the operator to be *bounded* if for any such $\alpha$ there is some (necessarily non-negative) real number $c$ such that

$$\|A\alpha\| \le c \|\alpha\| \tag{46}$$

Here, we are allowed to choose $c$ as large as we like, but it must be finite. Then we can define the *supremum (operator) norm* as the smallest possible choice of $c$ such that this expression (46) is valid for all $\alpha$ in the domain $D$. Equivalently, we can write[34]

$$\|A\|_{sup} = \sup_{\substack{\alpha \text{ in } D \\ \alpha \ne 0}} \frac{\|A\alpha\|}{\|\alpha\|} \tag{47}$$

With relation (46), we can also therefore write

$$\|A\alpha\| \le \|A\|_{sup} \|\alpha\| \tag{48}$$

since $\|A\|_{sup}$ is defined as the smallest possible value of $c$ for which Eq. (46) always works. In fact, a relation of this form Eq. (48) is a requirement for any operator norm, so quite generally for any operator norm $\|A\|$ we will require

$$\|A\alpha\| \le \|A\|\|\alpha\| \tag{49}$$

and such an expression (49) is useful in later proofs. Specifically, we will show this kind of relation also holds for the Hilbert-Schmidt norm that we introduce later.

Note here that the norm in $\|A\alpha\|$ is the vector norm as in Eq. (28) (though note formally here that this is the vector norm in the range $R$, so it would be based on the underlying inner product in the range space). In words, this is saying that this supremum norm for the operator A is the size of the "largest" possible vector we could produce in the range when starting with a unit length vector in the domain. By "largest" here, we mean the supremum (lowest possible upper bound) on the norm of the vector produced in the range.

Note that, with the definition of an operator norm, it becomes possible to consider the convergence not only of real numbers and vectors, but also of operators, and this will be important below.

## 6.2    Representing operators as matrices

Because any Hilbert space has some complete basis set, we can use this property and the underlying inner product to represent a linear operator as a matrix. This has important practical uses, of course, but it also has algebraic uses and helps us further define and extend the Dirac notation as being a particularly useful notation for Hilbert spaces and linear operators.

---

[34] In words, this notation means "the supremum of the number $\|A\alpha\| / \|\alpha\|$ for any possible choice of a non-zero vector or function $\alpha$ in the domain $D$ of the operator A".

Suppose, then, that we have two Hilbert spaces, $H_1$ and $H_2$; these may be the same Hilbert space, but here we also want to allow for the possibility that they are different. (We will need this when, for example, we are considering "source" and "receiving" spaces for waves.) We can propose vectors $\eta$ in $H_1$ and $\sigma$ and $\mu$ in $H_2$. We will also presume an orthonormal basis $\{\alpha_1, \alpha_2, \ldots\}$ in $H_1$ and an orthonormal basis $\{\beta_1, \beta_2, \ldots\}$ in $H_2$. Both Hilbert spaces may be infinite dimensional, and so these basis sets may also be infinite.

We presume that a bounded linear operator[35] $\mathsf{A}_{21}$ maps from vectors in space $H_1$ to vectors in space $H_2$, for example, mapping an vector $\eta$ in $H_1$ to some vector $\sigma$ in $H_2$

$$\sigma = \mathsf{A}_{21}\eta \tag{50}$$

Quite generally, we could construct the (underlying) inner product between this resulting vector and an arbitrary vector $\mu$ in $H_2$. Specifically, we would have

$$(\mu, \sigma)_2 \equiv (\mu, \mathsf{A}_{21}\eta)_2 \tag{51}$$

Note that this inner product is taken in $H_2$ (where we remember, as in Eq. (50), that $\mathsf{A}_{21}\eta$ is a vector in $H_2$), and we have used the subscript "2" to make this clear. This inner product is in the form of the underlying inner product in $H_2$.

Note again that the forms of the underlying inner products in the two different spaces $H_1$ and $H_2$ do not have to be the same; they just both have to be legal inner products. So the inner product in space $H_1$ might be a non-weighted inner product useful for representing, say, current sources, and that in $H_2$ might be a power or energy inner product for waves that could therefore be a weighted inner product. These possible differences in inner products in the two spaces mean that, for the moment, that we have to be careful to keep track of what space an inner product is in[36].

Now it will be useful to define what we will call a matrix element of the operator $\mathsf{A}_{21}$. In the most general situation which we are considering, where $H_1$ and $H_2$ could be different spaces, with different basis sets, we can define this *matrix element*, which is generally a complex number, as

$$a_{jk} = (\beta_j, \mathsf{A}_{21}\alpha_k)_2 \tag{52}$$

Again, this inner product is being taken in $H_2$, as indicated with the subscript "2".

Now let us consider an expression of the form Eq. (51) again, but this time we are going to represent each of the vectors $\eta$ and $\mu$ by expanding them on their corresponding basis sets using the underlying inner product in each space. So, we have

$$\eta = \sum_k r_k \alpha_k \tag{53}$$

and

---

[35] Note that the order of the subscripts on the operator $\mathsf{A}_{21}$ here is one that makes sense when we think of an operators operating on a function or vector on the "right", in space 1, to generate a new vector or function on the "left", in space 2 (which may be different from space 1). Indeed, for differential operators this "right to left" order is almost always implicit in the notation. Unless we invent a new notation, differential operators only operate to the right. Matrix operators can operate in either direction, but it is more conventional to think of column vectors as being the "usual" notation and row vectors as being an "adjoint" notation, in which case matrix-vector operations are also typically written in this same order.

[36] Note that it is generally meaningless to try to form an inner product between a function in one Hilbert space and a function in another Hilbert space; an inner product is a characteristic of a given Hilbert space, so we only need to put one subscript on the inner product in Eq. (51) to indicate the space in which it is being taken.

$$\mu = \sum_j t_j \beta_j \tag{54}$$

where the $r_k$ and the $t_j$ are complex numbers given by

$$r_k = (\alpha_k, \eta)_1 \tag{55}$$

(an inner product formed using the underlying inner product in $H_1$) and

$$t_j = (\beta_j, \mu)_2 \tag{56}$$

(an inner product formed using the underlying inner product in $H_2$). Then, we can rewrite Eq. (51) as

$$
\begin{aligned}
(\mu, A_{21}\eta)_2 &= \sum_j t_j^* \left( \beta_j, A_{21} \left[ \sum_k r_k \alpha_k \right] \right)_2 \\
&= \sum_{j,k} t_j^* r_k \left( \beta_j, A_{21}\alpha_k \right)_2 \\
&= \sum_{j,k} t_j^* a_{jk} r_k
\end{aligned}
\tag{57}
$$

Now we are in a position to make an "algebraic shift" towards a matrix-vector algebra, written in Dirac notation. Now we algebraically regard the "bra" vector $\langle \mu |$ as the row vector of expansion coefficients

$$\langle \mu | \equiv [t_1^*, t_2^*, \dots] \equiv \begin{bmatrix} t_1 \\ t_2 \\ \vdots \end{bmatrix}^{\dagger} \equiv \left( |\mu\rangle \right)^{\dagger} \tag{58}$$

which is equivalent to the "ket" version

$$|\mu\rangle \equiv \begin{bmatrix} t_1 \\ t_2 \\ \vdots \end{bmatrix} \tag{59}$$

and similarly the "ket" vector $|\eta\rangle$ is regarded as the column vector of expansion coefficients

$$|\eta\rangle \equiv \begin{bmatrix} r_1 \\ r_2 \\ \vdots \end{bmatrix} \tag{60}$$

Once we are working with these bra and ket vectors, we can also decide to regard an operator $A_{21}$ in algebraic expressions with bra and ket vectors as the matrix

$$A_{21} \equiv \begin{bmatrix} a_{11} & a_{12} & \cdots \\ a_{21} & a_{22} & \cdots \\ \vdots & \vdots & \ddots \end{bmatrix} \tag{61}$$

Then the sum $\sum_{j,k} t_j^* a_{21} r_k$ can be interpreted as the vector-matrix-vector product

$$\sum_{j,k} t_j^* a_{jk} r_k \equiv \langle \mu | A_{21} | \eta \rangle \tag{62}$$

Explicitly, we note that, quite generally, from Eqs. (57) and (62)

$$(\mu, A_{21}\eta)_2 = \langle \mu | A_{21} | \eta \rangle \tag{63}$$

The actual "underlying" operator $A_{21}$ operating on a function $\eta$ in $H_1$, as in the expression $A_{21}\eta$ inside the underlying inner product in $H_2$ on the left of Eq. (63), is only specified when it is operating "to the right"[37]; the expression "$\mu A_{21}$" does not necessarily have any meaning. However, once we have made this algebraic shift to the matrix-vector Dirac notation, the matrix-vector product $\langle \mu | A_{21}$ (which results in a row vector) is just as meaningful as the product $A_{21} | \eta \rangle$ (which results in a column vector). The fact that the underlying operator $A_{21}$ possibly only operates to the right has been "hidden" inside the matrix elements

$$a_{jk} \equiv \left( \beta_j, A_{21} \alpha_k \right)_2 \tag{64}$$

We could be criticized here for using the same notation for the matrix version of the operator and for the underlying linear operator, but there need be no confusion; if we see an expression such as $A_{21} \eta$, we are dealing with the underlying operator, which possibly only operates to the right, but if we see an expression such as $\langle \mu | A_{21}$, $A_{21} | \eta \rangle$, or $\langle \mu | A_{21} | \eta \rangle$, with the vectors in bra and/or ket form, then we are dealing with the matrix version of the operator.

In most use of Dirac notation, as, for example, in quantum mechanics, it is much more typical to have the operators map from a given Hilbert space to itself. Additionally, inner products other than a simple Cartesian form are unusual in quantum mechanics. Hence much of the subtlety we have been setting up here, in being careful about what inner product is in what space, and what form of inner product we are using, is unnecessary in quantum mechanics. Here, however, because we want to get the algebraic benefits of Dirac or matrix-vector algebra and we may well be operating between different Hilbert spaces with different inner products in each, we needed to set up this algebra with some care. The good news is that, with our understanding of how to use the underlying inner products in each space to evaluate expansion coefficients, as in Eqs. (55) and (56), and matrix elements, as in Eq. (64), we can make this algebraic shift to matrix-vector or Dirac notation and use their full power even in this more general situation..

We can usefully go one step further with Dirac notation here. We can also write the matrix $A_{21}$ itself in terms of bra and ket vectors. Again this is standard in other uses of Dirac notation, though, at least for the moment, we will be explicit about what spaces the vector are in by using "1" and "2" subscripts on the vectors. Specifically, we can write

$$A_{21} \equiv \sum_{j,k} a_{jk} \left| \beta_j \right\rangle_2 {}_1 \left\langle \alpha_k \right| \tag{65}$$

Then

$$
\begin{aligned}
{}_2\langle \mu | A_{21} | \eta \rangle_1 &\equiv {}_2 \left\langle \mu \left| \left( \sum_{j,k} a_{jk} \left| \beta_j \right\rangle_2 {}_1 \left\langle \alpha_k \right| \right) \right| \eta \right\rangle_1 \\
&= \sum_p t_p^* {}_2 \left\langle \beta_p \left| \left( \sum_{j,k} a_{jk} \left| \beta_j \right\rangle_2 {}_1 \left\langle \alpha_k \right| \right) \sum_q r_q \left| \alpha_q \right\rangle_1 \right. \\
&= \sum_p t_p^* \sum_{j,k} \delta_{pj} a_{jk} {}_1 \delta_{kq} \sum_q r_q \\
&= \sum_{j,k} t_j^* a_{jk} r_k
\end{aligned}
\tag{66}
$$

which is again the same as the result in the original equation Eq. (57), so this approach for writing matrices works here also.

Quite generally, a form like $\left| \beta_j \right\rangle_2 {}_1 \left\langle \alpha_k \right|$ is an *outer product*. In contrast to the inner product, which produces a complex number from the multiplication in "row vector – column vector" order, and which necessarily only involves vectors in the same Hilbert space, the outer product can be regarded as generating a matrix from the multiplication in "column vector – row vector" order, and can involve vectors in different Hilbert spaces. Dropping the additional subscript notation on the vectors, instead of Eq. (65) we will just write

$$A_{21} \equiv \sum_{j,k} a_{jk} \left| \beta_j \right\rangle \left\langle \alpha_k \right| \tag{67}$$

---

[37] For example, derivative operators are usually only defined as operating to the right.

Note that a linear operator like $A_{21}$ from one Hilbert space to another can be written in such an outer product form as in Eq. (65) on any desired basis sets for each Hilbert space. Of course, the numbers $a_{jk}$ will be different depending on the basis sets chosen.

This statement of the operator as a matrix in Dirac notation completes our "algebraic shift". From this point on, we use either the notation with functions written as just Greek letters such as $\alpha$ and $\beta$ with (underlying) inner products $(\alpha, \beta)$, or Dirac notation with functions written as kets, such as $|\alpha\rangle$ and $|\beta\rangle$ (or their corresponding bra versions $\langle\alpha|$ and $\langle\beta|$) and (Cartesian) inner products written as $\langle\alpha|\beta\rangle \equiv \langle\alpha\|\beta\rangle$. Importantly, because the underlying inner products are always used in constructing the vectors and matrices in the Dirac notation, the result of any such expression in both notations is the same. So, we can move between notations depending on convenience, and we will do so below. The ability to use the associative property of matrix-vector notation (including "breaking up" the inner product as in $\langle\alpha|\beta\rangle \equiv \langle\alpha\|\beta\rangle$) often results in considerable algebraic simplification.

## 6.3 Adjoint operator

We can now usefully define what we will call the *adjoint operator* of an operator $A_{21}$. We will (temporarily) call this operator $B_{12}$ so we can define its properties without some confusion of notation[38].

We define this adjoint operator through the relation, for any vectors $\eta$ in $H_1$ and $\mu$ in $H_2$,

$$(\mu, A_{21}\eta)_2 = (B_{12}\mu, \eta)_1 \tag{68}$$

Note that $B_{12}$ is an operator that maps from Hilbert space $H_2$ to Hilbert space $H_1$. Note, too, that in this case, the inner product on the right hand side is performed in $H_1$; both $B_{12}\mu$ and $\eta$ are vectors in $H_1$. Now, similarly to Eq. (52), we will write a "matrix element" between the appropriate basis functions, called for the moment

$$b_{kj} = (\alpha_k, B_{12}\beta_j)_1 \tag{69}$$

Now that we have these matrix elements for $B_{12}$ defined, we can make the algebraic shift to matrix-vector algebra. We treat $B_{12}$ as a matrix with matrix elements as in Eq. (69) and we write the vectors of expansion coefficients for $\mu$ and $\eta$ as in Eqs. (59) and (60), respectively. So, instead of Eq. (68) we can write

$$\langle\mu|A_{21}|\eta\rangle = \left(B_{12}|\mu\rangle\right)^\dagger|\eta\rangle = \langle\mu|B_{12}^\dagger|\eta\rangle \tag{70}$$

where the "$\dagger$" is the matrix and vector Hermitian adjoint operation, as discussed in section 4.3, and where we used the known standard result for matrix-vector multiplication that

$$\left(C|\theta\rangle\right)^\dagger \equiv \left(|\theta\rangle\right)^\dagger C^\dagger = \langle\theta|C \tag{71}$$

for a matrix $C$ and a vector $|\theta\rangle$. Now expanding the vectors on their basis sets on the right hand side of Eq. (70), we have

$$\langle\mu|A_{21}|\eta\rangle = \sum_j t_j^* \left\langle\beta_j\left|B_{12}^\dagger\left(\sum_k r_k|\alpha_k\rangle\right)\right.\right\rangle$$
$$= \sum_{j,k} t_j^* r_k \left\langle\beta_j\left|B_{12}^\dagger\right|\alpha_k\right\rangle \tag{72}$$

Now, if the matrix $B_{12}$ has matrix elements $b_{jk}$ in the $j$th row and $k$th column, then the matrix $B_{12}^\dagger$ has matrix elements $b_{kj}^*$ in the $j$th row and $k$th column, i.e.,

$$\left\langle\beta_j\left|B_{12}^\dagger\right|\alpha_k\right\rangle = b_{kj}^* \tag{73}$$

So from Eq. (72)

---

[38] A confusion that could make it seem that we are assuming what we are trying to prove

$$\langle \mu | A_{21} | \eta \rangle = \sum_{j,k} t_j^* b_{kj}^* r_k \tag{74}$$

From Eq. (66), the left hand side of Eq. (74) has to equal $\sum_{j,k} t_j^* a_{jk} r_k$. Hence we have

$$\sum_{j,k} t_j^* a_{jk} r_k = \sum_{j,k} t_j^* b_{kj}^* r_k \tag{75}$$

However, the vectors or functions $\eta$ and $\mu$ (and hence also the sets of coefficients $r_k$ and $t_j$) are arbitrary in their Hilbert spaces. So therefore we must have

$$b_{kj}^* = a_{jk} \tag{76}$$

which means that the adjoint operator $B_{12}$ is (at least in matrix form), the Hermitian adjoint of the original operator $A_{21}$.

$$B_{12} \equiv A_{21}^\dagger \tag{77}$$

so we can write as a defining equation of an adjoint operator

$$(\mu, A\eta) = (A^\dagger \mu, \eta) \tag{78}$$

for any vectors $\eta$ and $\mu$ in the appropriate Hilbert spaces. (Here we have dropped the subscripts for simplicity of notation.) Note that this expression Eq. (78) can be stated for the general case of the operator $A$, not just its matrix representation. Note also that we can see from this matrix form that

$$\left(A^\dagger\right)^\dagger = A \tag{79}$$

So, henceforth, we can write the adjoint operator to $A_{21}$ as simply $A_{21}^\dagger$, and our adjoint operator is simply the Hermitian adjoint of the original operator[39]. Note that we have proved this even for different spaces $H_1$ and $H_2$ with possibly different inner products in both spaces.

## 6.4   Identity operator

The mathematics of expanding on an orthonormal basis and making the algebraic shift to matrix-vector (i.e., Dirac) notation gives us an algebraically useful form for the identity operator in a space. To see this, we start by making an algebraic shift for the basis vectors, writing the basis functions $\alpha_j$ in a given space themselves as Dirac bra $|\alpha_j\rangle$ or ket $\langle \alpha_j|$ vectors[40]. Hence, in Dirac notation, we can write for some function $\gamma$ in the same Hilbert space

$$|\gamma\rangle = \sum_j \langle \alpha_j | \gamma \rangle | \alpha_j \rangle \tag{80}$$

Now $\langle \alpha_j | \gamma \rangle$ is just a complex number, so we can move it within the expression in the sum to obtain

$$|\gamma\rangle = \sum_j |\alpha_j\rangle \langle \alpha_j | \gamma \rangle = \left( \sum_j |\alpha_j\rangle \langle \alpha_j | \right) |\gamma\rangle \tag{81}$$

where now we have explicitly split up the inner product into a product of a "bra" and a "ket" vector. Using the associative properties of matrix-vector multiplications, inserting the parentheses in the expression on the far right of Eq. (81), we now have the outer product $\left|\alpha_j\right\rangle\left\langle\alpha_j\right|$ appearing in the sum. In this case, we see that the effect of the operator

$$\mathsf{I}_{op} = \sum_j \left|\alpha_j\right\rangle\left\langle\alpha_j\right| \tag{82}$$

is that it acts as *identity operator* for all vectors $\gamma$ in this space. Note that the identity operator can be written as such a sum of such outer products[41] using any (complete) basis set in the space[42], a property that is mathematically algebraically very useful in proofs and other manipulations.

---

[39] Note that, though this adjoint operator $A_{21}^\dagger$ is written with the subscripts in the order, from left to right, "2-1", it is an operator that maps from $H_2$ to $H_1$; changing the order here would have created possibly more confusion.

[40] We can regard the basis functions themselves as being expanded on a basis (possibly, but not necessarily, a different basis), using the underlying inner product to calculate the expansion coefficients, just as for any other function.

## 6.5    Compact operators

Compact operators are a category of operators that go beyond just being bounded. The notion of compactness is essentially one we need when working with operators in infinite dimensional spaces[43]. Now[44], we will build towards an understanding of what a compact operator is, which in turn can give us a sense of their importance in physical problems[45].

### 6.5.1    Definition of a compact operator

The formal definition of a *compact (linear) operator* is somewhat abstract.

> We presume we have a linear operator $\mathsf{A}$ that maps from one normed (vector) space, *F,* to another, possibly different vector space *G*. So a vector $\alpha$ in *F* leads to a corresponding vector $\gamma = \mathsf{A}\alpha$ in *G*. A compact (linear) operator is such an operator $\mathsf{A}$ with the additional property that, for any set of bounded vectors $\{\alpha_m\}$ in *F*, the set of corresponding vectors $\gamma_m$ in *G* is precompact (i.e., one whose closure is compact). (83)

From this definition we can prove a theorem that gives a somewhat more direct criterion for an operator to be compact. Specifically

> The operator $\mathsf{A}$ (from the normed space *F* to the normed space *G*) is compact if and only if it maps every bounded sequence $(\alpha_m)$ of vectors in *F* into a sequence in *G* that has a convergent subsequence. (84)

See 11.4 "Proof (4) of a criterion for compactness of an operator" below.

### 6.5.2    An illustrative extreme example

To see the power of, and the need for, operator compactness, consider the following "extreme" example. Consider an infinite-dimensional Hilbert space, with an orthonormal basis $\{\alpha_1, \alpha_2, \ldots\}$. For any two such basis vectors, the "distance" between them, as defined by the metric (see section 4.5), is

$$d(\alpha_j, \alpha_k) \equiv \sqrt{(\alpha_j - \alpha_k, \alpha_j - \alpha_k)} = \sqrt{(\alpha_j, \alpha_j) + (\alpha_k, \alpha_k) - (\alpha_k, \alpha_j) - (\alpha_j, \alpha_k)}$$
$$= \sqrt{1 + 1 - 0 - 0} = \sqrt{2}$$

(85)

---

[41] Note therefore that the identity operator $\mathsf{I}_{op}$ can then be considered as a sum of these "outer product" matrices.

[42] In general, different spaces have different identity operators, and so, if necessary, we can subscript the identity operator to indicate what space it operates in.

[43] Compactness is a somewhat trivial property in bounded finite dimensional spaces because all bounded finite dimensional linear operators are compact, as we will prove below.

[44] The mathematical definitions, theorems, and proofs on compact operators in this section are based on Kreyszig's approach [2], especially theorems 8.1-5, 8.1-4 (a), 2.5-3 (in particular, the proof of the compactness of any closed and bounded finite dimensional normed space), and 2.4-1, though we have harmonized the notation with our approach and avoided introducing some concepts that are not required elsewhere in our discussion.

[45] In the physics of waves, the properties of compact operators are behind the notion of diffraction limits and limitations on the number of usable channels in communications, for example.

(This can be visualized as the distance between the "tips" of two unit vectors that are at right angles.) So, we can construct an infinite sequence that is just the basis vectors, each used exactly once, such as the sequence $(\alpha_1, \alpha_2, \ldots)$. This sequence does not converge, and has no convergent subsequences[46]; every pair of elements in the sequence has a "distance" between them of $\sqrt{2}$. A compact operator operating on that infinite sequence of different basis vectors will get rid of this problem in the vectors it generates – those will have some convergent subsequence. So, the compact operator eliminates one troubling aspect of working with infinite dimensional spaces.

### 6.5.3    Compactness of operators with a finite-dimensional range

As we proceed to understand compact operators, first, we note that any operator with a finite dimensional range is compact[47]. The proof of this statement itself has several steps to it, and is given below in 11.5 "Proof (5) of compactness of operators with finite dimensional range". We then use this first result to prove the following theorem[48]:

> Consider an infinitely long sequence $(A_n)$ of compact linear operators from a normed space $F$ (which can, of course, be a Hilbert space) into a Hilbert space $G$. If $\|A_n - A\| \to 0$ as $n \to \infty$ for some operator $A$, then this limit operator $A$ is compact.

(86)

We give the proof of this theorem below in 11.6 "Proof (6) of compact limit operator for convergent compact operators". This theorem (86) is a statement that, if an operator can be approximated to any degree of accuracy by a finite matrix, then it is compact. This is not a statement that any compact operator can be approximated by a sufficiently large matrix, but the broad class of Hilbert-Schmidt (compact) operators can be rigorously approximated in this way, and we develop the theory of them below.

## 6.6    Hilbert-Schmidt operators

It is possible to derive various further results for compact operators, and we will do continue to do so. However, for a broad class of physical problems[49], we will be particularly interested in so-called Hilbert-Schmidt operators; the mathematics of these is somewhat simpler than carrying forward the full general mathematics of compact operators. They are perhaps also easier to understand intuitively, so we restrict some of our discussion to them. We first have to define what we mean by a Hilbert-Schmidt operator, and then we will prove that all Hilbert-Schmidt operators are compact.

### 6.6.1    Definition of a Hilbert-Schmidt operator and the sum-rule limit

We can define a *Hilbert-Schmidt operator* as follows.

> We presume we have a Hilbert space $H_1$ with an orthogonal basis $\{\alpha_1, \alpha_2, \ldots\}$ (which we presume to be orthonormal for convenience), and a bounded operator $A$ that maps from vectors in $H_1$ to vectors in a possibly different Hilbert space[50] $H_2$. Then $\mathsf{A}$ is a Hilbert-Schmidt operator if and only if

$$S \equiv \sum_j \left\| \mathsf{A}\alpha_j \right\|^2 < \infty \tag{87}$$

---

[46] The same problem does not arise in finite-dimensional spaces; if we construct an infinitely long sequence made up from just the finite number of basis vectors in the space, we will have to repeat at least one of the basis vectors an infinite number of times, which gives us at least one convergent subsequence – the sequence consisting of just that basis vector repeated an infinite number of times.

[47] The reader may already be able to see this informally and intuitively from the above "extreme" example and the preceding footnote[46].

[48] This theorem is a somewhat restated version of Theorem 8.1.5 in Kreyszig [2], and we give a version of that proof.

[49] For example, essentially all the "Green's function" operators we encounter in dealing with the physics of waves generated by sources correspond to Hilbert-Schmidt operators.

Since the result of this sum is finite, we can give it a name and a notation, calling it[51] the *sum rule limit S*, subscripted if necessary to show it is associated with some specific operator. The square root of this (necessarily non-negative) sum-rule limit $S$ can be called the *Hilbert-Schmidt norm* of the operator, i.e.,

$$\|\mathsf{A}\|_{HS} = \sqrt{S} \equiv \sqrt{\sum_j \left\| \mathsf{A}\alpha_j \right\|^2} \tag{88}$$

For any arbitrary complete basis sets $\left\{ \left| \alpha_j \right\rangle \right\}$ and $\left\{ \left| \beta_k \right\rangle \right\}$ in $H_1$, starting from this definition, we can prove three other equivalent expressions[52] for $S$, given in the three lines in the equations (89) below

$$
\begin{aligned}
S \equiv \|\mathsf{A}\|_{HS}^2 &= \sum_j \left\langle \alpha_j \left| \mathsf{A}^\dagger \mathsf{A} \right| \alpha_j \right\rangle = \sum_k \left\langle \beta_k \left| \mathsf{A}^\dagger \mathsf{A} \right| \beta_k \right\rangle \\
&= \sum_{j,k} \left| a_{kj} \right|^2 \\
&\equiv Tr\left( \mathsf{A}^\dagger \mathsf{A} \right) = Tr\left( \mathsf{A}\mathsf{A}^\dagger \right)
\end{aligned}
\tag{89}
$$

See 11.7 "Proof (7) of equivalent statements of the Hilbert-Schmidt sum rule limit S". Since all of these different statements of $S$ are equivalent, proving that any one of these versions is finite on any complete basis is sufficient to prove an operator $\mathsf{A}$ is a Hilbert-Schmidt operator. We can also now explicitly prove that the required property for any operator norm as given in relation (49) ($\|\mathsf{A}\alpha\| \le \|\mathsf{A}\|\|\alpha\|$) also applies for this Hilbert-Schmidt norm. See 11.8 "Proof (8) of the operator norm inequality for the Hilbert-Schmidt norm".

### 6.6.2 Compactness of Hilbert-Schmidt operators

One particularly important property of Hilbert-Schmidt operators is that they are compact. The proof of this is given below in 11.9 "Proof (9) of compactness of Hilbert-Schmidt operators".

### 6.6.3 Approximating Hilbert-Schmidt operators by sufficiently large matrices

Consider a Hilbert-Schmidt operator $\mathsf{A}$ that maps from a Hilbert space $H_1$ to a possibly different Hilbert space $H_2$. Consider also an $m \times n$ matrix $\mathsf{A}_{mn}$ that is a "truncated" version of the matrix version of $\mathsf{A}$.

---

[50] These Hilbert spaces can be infinite dimensional.

[51] This "sum rule limit" name is one we are creating, and it not standard in the mathematics literature.

[52] The Hilbert-Schmidt norm is often also written in integral form. Indeed, once we consider physical operators like Green's functions, this is very appropriate. Here, for the purposes of our mathematics we mostly omit that, regarding it as a special case of forms derived from the infinite sum as in Eq. (89). If we write it out as an integral, we have to be more specific about the form of the corresponding operator, such as a Green's function that might be operating on different kinds of physical spaces (e.g., 1-dimensional or 3-dimensional), and it might have some more sophisticated character, including tensor or dyadic forms. For completeness, though, one specific example, for a scalar Green's function $G\left(\mathbf{r}_2;\mathbf{r}_1\right)$ giving the scalar wave a position $\mathbf{r}_2$ in volume $V_2$ in response to a point source at position $\mathbf{r}_1$ in volume $V_1$, would be $S = \int_{V_2}\int_{V_1} \left| G\left(\mathbf{r}_2;\mathbf{r}_1\right) \right|^2 d^3\mathbf{r}_1 d^3\mathbf{r}_2$. See [1] for more discussion of such physical Green's functions. Indeed, whether a specific operator is a Hilbert-Schmidt one will often be determined by such an integral. An important point is that, as a result, a very broad class of Green's function operators, including those in wave problems, are Hilbert-Schmidt operators. To justify that more fully, we need to consider the physics behind such operators; situations with finite volumes, and where the response from a finite source is itself finite, are, however, generally going to correspond to Hilbert-Schmidt operators [1]. It is that finiteness from the physics that allows us to exploit the mathematics of compact operators, and especially Hilbert-Schmidt ones.

Then we can prove that the vector result $A_{mn}\mu$ of operating with $A_{mn}$ on any vector $\mu$ in $H_1$ converges to the vector result $A\mu$ if we take $m$ and $n$ to be sufficiently large. See the 11.10 "Proof (10) of approximation of Hilbert-Schmidt operators by sufficiently large matrices" below. Hence, Hilbert-Schmidt operators can always be approximated by sufficiently large finite matrices.

### 6.6.4 Nature of operators $A^\dagger$, $A^\dagger A$ and $A A^\dagger$ for a Hilbert-Schmidt operator A

For a Hilbert-Schmidt operator $A$, the operators $A^\dagger$, $A^\dagger A$ and $AA^\dagger$ are Hilbert-Schmidt operators (and are therefore also compact).

(90)

See 11.11 "Proof (11) of Hilbert-Schmidt and compact nature of various operators derived from Hilbert-Schmidt operators".

## 6.7 Hermitian operators

A particularly powerful and useful class of operators are the so-called "Hermitian" operators. The most general definition of a *Hermitian* or *self-adjoint* operator $A$ is that, for all vectors or functions $\beta$ and $\gamma$ in the relevant Hilbert space or spaces,

$$(\beta, A\gamma) = (A\beta, \gamma)$$ (91)

If we compare this with the definition of the adjoint operator, Eq. (68), we see that this means this operator is equal to its own adjoint. Equivalently then, in particular if we are considering the matrix representation of the operator on some basis,

$$A = A^\dagger$$ (92)

and for the matrix elements of the operator

$$a_{jk} = a_{kj}^*$$ (93)

An equivalent statement would therefore be that this matrix is equal to its own "conjugate transpose".

### 6.7.1 Hermiticity of $A^\dagger A$ and $A A^\dagger$

In general, the operators that map from, say, sources in one space to waves in another are often not Hermitian, so we will certainly be dealing with non-Hermitian operators. However, the operators $A^\dagger A$ and $AA^\dagger$ are Hermitian, as is simply proved using the usual rule for the Hermitian adjoint of a product being the "flipped round" product of the Hermitian adjoints, and the fact that the Hermitian adjoint of a Hermitian adjoint takes us back to the original matrix, i.e.,

$$\left(A^\dagger A\right)^\dagger = \left(A\right)^\dagger \left(A^\dagger\right)^\dagger = A^\dagger A$$ (94)

$$\left(AA^\dagger\right)^\dagger = \left(A^\dagger\right)^\dagger \left(A\right)^\dagger = AA^\dagger$$ (95)

Hence:

For a Hilbert-Schmidt operator $A$ (which is not necessarily Hermitian), the operators $A^\dagger A$ and $AA^\dagger$ are both compact and Hermitian

(96)

This result is very important in the mathematics here, as will become apparent once we understand the properties of eigenvalues and eigenfunctions of compact Hermitian operators.

# 7 Eigenvalues and eigenfunctions of compact Hermitian operators

Here we will develop the theory of eigenvalues and eigenfunctions of compact Hermitian operators. There are several important ultimate results. One is that the eigenfunctions are all orthogonal and form a complete set (or one that can be easily completed). Another is that the eigenfunctions are progressively the ones that maximize the possible "strengths" of the effects of the operator on functions. The detail of this will become clear as we develop this, leading up to the so-called "spectral theorem" that gives the final results we want. There are several steps to this, and we need to start with some preliminary results.

As we go through some of these operator properties, including some "eigen properties", we will label some of them as (OE1), (OE2), etc.

## 7.1 A property of Hermitian operators

First, we prove a property of a Hermitian operator $A$ . By definition, as in Eq. (91), for any vector $\beta$ ,

$$(\beta, A\beta) = (A\beta, \beta) \tag{97}$$

But the property (IP3) of an inner product requires that

$$(\beta, A\beta) = (A\beta, \beta)^* \tag{98}$$

So, $(A\beta, \beta) = (A\beta, \beta)^*$ , which therefore requires that $(A\beta, \beta)$ is real, and hence also $(\beta, A\beta)$ is real. So, quite generally,

(OE1)        for a Hermitian operator $(\beta, A\beta)$ is a real number    (99)

## 7.2 Definition of eigenfunctions and eigenvalues

Quite generally, for an operator $A$ , some vector $\alpha$ is an *eigenvector* (or *eigenfunction*) of $A$ if and only if

(OE2)                    $A\alpha = c\alpha$    (100)

where $c$ is some number, possibly complex, that is then called the *eigenvalue* associated with this eigenvector.

## 7.3 Reality of eigenvalues of a Hermitian operator

For any eigenvector $\alpha$ ,

$$(\alpha, A\alpha) = (\alpha, c\alpha) = c(\alpha, \alpha) \tag{101}$$

where we used property (IP2) of an inner product. But $(\alpha, \alpha)$ is necessarily real by inner product property (IP3), and from (99) (OE1), $(\beta, A\beta)$ is real for any vector $\beta$ , so therefore the eigenvalue $c$ is also real; so,

(OE3)        All eigenvalues of Hermitian operators are necessarily real.    (102)

## 7.4 Finite number of eigenvectors of a compact operator for a given non-zero eigenvalue

It is possible for a given eigenvalue to have multiple different linearly independent eigenvectors. The dimensionality of the space that includes all such eigenvectors for a given eigenvalue can be called the *multiplicity* or the *degeneracy* of the eigenvalue. We can state an important property:

(OE4) A non-zero eigenvalue of a compact Hermitian operator has finite multiplicity (103)

We prove this below in 11.13 "Proof (13) of finite multiplicity". In 11.14 "Proof (14) that the eigenvalues of Hermitian operator on an infinite dimensional space tend to zero", we also show that

(OE5)   If a compact Hermitian operator is operating on an infinite dimensional space, then
the sequence of eigenvalues $(c_p)$ must tend to zero as $p \to \infty$. (104)

## 7.5    Orthogonality of eigenvectors for different eigenvalues

For a compact Hermitian operator $A$, suppose it has two different eigenvalues $a$ and $b$. So for some corresponding (non-zero) eigenvectors $\alpha$ and $\beta$, we have

$$A\alpha = a\alpha$$

and

$$A\beta = b\beta$$

Then

$$b(\alpha, \beta) = (\alpha, b\beta) = (\alpha, A\beta) \tag{105}$$

Now, by the Hermiticity of $A$, we have

$$(\alpha, A\beta) = (A\alpha, \beta) = (a\alpha, \beta) = a^*(\alpha, \beta)$$
$$= a(\alpha, \beta) \tag{106}$$

where in the last step we used the fact that the eigenvalues of Hermitian operators are real, as in (102) (OE3) above. Hence, putting the left hand side of Eq. (105) together with the right hand side of Eq. (106) and rearranging, we have

$$(b - a)(\alpha, \beta) = 0 \tag{107}$$

Since by assumption $b \neq a$, then $(\alpha, \beta) = 0$, which means

(OE6)  for a Hermitian operator, eigenvectors for different eigenvalues are orthogonal (108)

## 7.6    A preliminary result for the supremum norm of compact Hermitian operators

So far, we have defined the supremum norm of an operator $A$ on a Hilbert space $H$ in terms of Eq. (47). For algebraic convenience, but without loss of generality, we rewrite this for vectors $\alpha$ of unit norm. So, we have

$$\|A\| = \sup_{\|\alpha\|=1} \|A\alpha\| \tag{109}$$

If $A$ is a compact Hermitian operator, we can prove that the supremum norm of $A$ can be rewritten as

$$\|A\| = \sup_{\|\alpha\|=1} |(\alpha, A\alpha)| \tag{110}$$

We prove[53] this below in 11.15 "Proof (15) of Hermitian operator supremum norm".

---

[53] For a similar proof, see [3], pp.198 − 199, Lemma 8.26. Our proof is not identical because we avoided requiring some prior results used in that proof, proving some parts directly instead, and we avoided some re-use of notation.

This result, Eq. (110), is at the core of the main results we will prove for eigenvectors of Hermitian operators. Note that, with the vector $\alpha$ also appearing on the left-hand side of the inner product $(\alpha, A\alpha)$, this result is saying, effectively, that the "largest" possible vector that can be produced by an operator acting on a unit-length vector is one that lies in the same or the opposite "direction" compared to the original vector, for some choice of that vector.

## 7.7 Spectral theorem for compact Hermitian operators

The spectral theorem essentially allows us to conclude that the eigenfunctions of a compact Hermitian operator form a complete set (or one that is easily completed), which is an extremely useful result. Along the way in this proof, we also establish an important maximization (or minimization) property of eigenvalues and their associated eigenvectors. The *spectral theorem* for the eigenfunctions of a compact Hermitian operator can be stated as follows:

> For a compact Hermitian operator $A$ mapping from a Hilbert space $H$ onto itself, the set of eigenfunctions $\{\beta_j\}$ of $A$ is complete for describing any vector $\phi$ that can be generated by the action of the operator on an arbitrary vector $\psi$ in the space $H$, i.e., any vector $\phi = A\psi$. If all the eigenvalues of $A$ are non-zero, then the set $\{\alpha_j\}$ will be complete for the Hilbert space $H$; if not, then we can extend the set by Gram-Schmidt orthogonalization to form a complete set for $H$.

(111)

(See the 11.3 "Proof (3) of the existence of a basis for a Hilbert space" for a discussion of Gram-Schmidt orthogonalization). We prove this theorem below in 11.16 "Proof (16) of the spectral theorem".

A consequence is that we can write any such compact Hermitian operator in terms of its eigenfunctions and corresponding eigenvalues as

$$A = \sum_{j=1}^{\infty} r_j \beta_j \left( \beta_j, \cdot \right)$$

(112)

(where we will substitute the vector being operated on for the dot "$\cdot$" when we use the operator) or, in Dirac notation

$$A = \sum_{j=1}^{\infty} r_j \left| \beta_j \right\rangle \left\langle \beta_j \right|$$

(113)

Here, the eigenvalues $r_j$ are whatever ones are associated with the corresponding eigenvector $\beta_j$. (Note in both Eqs (112) and (113) that, for the case of degenerate eigenvalues, we presume that we have written an orthogonal set of eigenvectors for each such degenerate eigenvalue (which we are always free to do) and for indexing purposes for an $p$-fold degenerate eigenvalue, we simply repeat it $p$ times in this sum, once for each of the corresponding eigenvectors.)

## 7.8 Maximization property of eigenvectors

As part of the proof of the spectral theorem (11.16 "Proof (16) of the spectral theorem"), we also show the following.

> The eigenvectors $\beta_j$ of a compact Hermitian operator can be found by a progressive variational technique, finding the largest possible result for $\left\| A\beta_j \right\|$ where $\beta_j$ is constrained to be orthogonal to all the previous eigenvectors. This will also give a corresponding set of eigenvalues $r_j$ in descending order of their magnitude.

(114)

This means, physically, that the eigenfunctions are essentially the "best" functions we can choose if we are trying to maximize performance in specific ways (such as maximizing power coupling between sources and the resulting waves), and we could even find them physically just by looking for the best such performance.

# 8  Inner products based on positive Hermitian operators

As we mentioned above in section 4.4, there is an even more general class of entities that have all the properties of inner products, and hence are inner products, and we can introduce those now. We need these ideas for defining "energy" inner products for electromagneticsm, for example [1].

## 8.1  Operator-weighted inner product

Suppose we have a (linear) Hermitian operator $\mathsf{A}$ that acts on functions such as $\alpha$, $\beta$, and $\gamma$, and suppose we already have defined an inner product of the form $(\beta, \alpha)$ with all the properties IP1 to IP4 as in (23). Now, the action of $\mathsf{A}$ on a vector $\gamma$ is to generate another vector $\alpha$ as in

$$\alpha = \mathsf{A}\gamma \tag{115}$$

So, we can form the inner product $(\beta, \alpha) \equiv (\beta, \mathsf{A}\gamma)$. Now, from the Hermiticity of $\mathsf{A}$, we know that $(\beta, \mathsf{A}\gamma) = (\mathsf{A}\beta, \gamma)$, as in Eq. (91), and by (IP3), we know that $(\beta, \mathsf{A}\gamma) = (\mathsf{A}\gamma, \beta)^*$. So, let us define a new entity, which we could call[54] an *operator-weighted inner product*,

$$(\beta, \gamma)_{\mathsf{A}} \equiv (\beta, \mathsf{A}\gamma) \tag{116}$$

Then, using first Eq. (91) and then IP3, we have

$$(\beta, \gamma)_{\mathsf{A}} = (\beta, \mathsf{A}\gamma) = (\mathsf{A}\beta, \gamma) = (\gamma, \mathsf{A}\beta)^* \equiv (\gamma, \beta)_{\mathsf{A}}^* \tag{117}$$

Hence this new entity, based on a Hermitian operator $\mathsf{A}$, also satisfies the property IP3 of an inner product. It is straight forward to show that, because $\mathsf{A}$ is linear, this entity also satisfies (IP1), as in

$$(\gamma, \alpha + \beta)_{\mathsf{A}} \equiv (\gamma, \mathsf{A}(\alpha + \beta)) = (\gamma, \mathsf{A}\alpha + \mathsf{A}\beta) = (\gamma, \mathsf{A}\alpha) + (\gamma, \mathsf{A}\beta)$$
$$= (\gamma, \alpha)_{\mathsf{A}} + (\gamma, \beta)_{\mathsf{A}} \tag{118}$$

and (IP2), as in

$$(\gamma, a\alpha)_{\mathsf{A}} \equiv (\gamma, \mathsf{A}a\alpha) = (\gamma, a\mathsf{A}\alpha) = a(\gamma, \mathsf{A}\alpha) \equiv a(\gamma, \alpha)_{\mathsf{A}} \tag{119}$$

As for (IP4), we already know that any entity $(\beta, \mathsf{A}\beta)$ is a real number, as shown in property (OE1) (Eq. (99)). However, it is not in general true for a Hermitian operator $\mathsf{A}$ that $(\beta, \mathsf{A}\beta)$ is positive. So for $(\beta, \gamma)_{\mathsf{A}}$ to be an inner product, we need one further restriction on $\mathsf{A}$, which is that it should be a *positive operator*, which by definition[55] means that

$$(\beta, \mathsf{A}\beta) \geq 0 \tag{120}$$

and hence by definition $(\beta, \gamma)_{\mathsf{A}}$ satisfies (IP4) with this restriction.

---

[54] As an explicit name, this "operator-weighted inner product" is a term we are creating here as far as we know, though this idea is known and this name may therefore be implicitly obvious.

[55] Note that there is some variation in notation in mathematics texts. Kreyszig [2] uses this definition for a positive operator, for example, and if the "$\geq$" sign is replaced by a "$>$" sign in (120), he would then call the operator *positive-definite*. Others, however, such as [5], would give (120) as the definition for a *non-negative operator*, using "positive operator" only if the "$\geq$" sign is replaced by a "$>$" sign.

So, for any positive (linear) Hermitian operator $\mathsf{A}$, we can construct an (operator-weighted) inner product of the form given by Eq. (116). (See also [4], p. 168.) The weighted inner product as in Eq. (26) can be viewed as a special case of this more general inner product[56].

## 8.2    Transformed inner product

One particular kind of positive operator is one that can be written in the form

$$\mathsf{A} = \mathsf{B}^{\dagger}\mathsf{B} \tag{121}$$

where $\mathsf{B}$ is a linear operator. We can prove that $\mathsf{A}$ in this case is a positive operator. Consider an arbitrary vector $\beta$, and form the inner product

$$\left(\beta, \mathsf{A}\beta\right) = \left(\beta, \mathsf{B}^{\dagger}\mathsf{B}\beta\right) = \left(\beta, \mathsf{B}^{\dagger}\left(\mathsf{B}\beta\right)\right) \tag{122}$$

But, by the defining property of an adjoint operator, as in Eq. (78), and with the property (79)

$$\left(\beta, \mathsf{B}^{\dagger}\left(\mathsf{B}\beta\right)\right) = \left(\mathsf{B}\beta, \mathsf{B}\beta\right) \tag{123}$$

which is the inner product of a vector with itself, which is necessarily greater than or equal to zero. So for an operator as in Eq. (121)

$$\left(\beta, \mathsf{A}\beta\right) \geq 0 \tag{124}$$

hence proving $\mathsf{A} = \mathsf{B}^{\dagger}\mathsf{B}$ is a positive operator. Hence, for any such operator, we could form an operator-weighted inner product.

We can, however, take an additional step that opens another sub-class of inner products. Specifically, we could define what we could call[57] a *transformed inner product*. We can regard the operator $\mathsf{B}$ as transforming[58] the vector $\beta$ - after all, $\mathsf{B}$ operating on $\beta$ is just a linear transform acting on $\beta$ – and we could write generally

$$\left(\beta, \gamma\right)_{T\mathsf{B}} \equiv \left(\mathsf{B}\beta, \mathsf{B}\beta\right) \tag{125}$$

where our subscript notation "$T\mathsf{B}$" indicates this inner product with respect to the transformation $\mathsf{B}$ of the vectors in the inner product. Our proof above, Eqs. (122) to (124) shows that this inner product $\left(\beta, \gamma\right)_{T\mathsf{B}}$ also satisfies (IP4).

# 9    Singular-value decomposition

The idea of singular-value decomposition (SVD), especially for finite matrices, is a well-known mathematical technique for rewriting matrices. As a general approach to rewriting linear operators, it may be less well known, but in wave problems [1][6][7] this approach can be particularly useful and physically meaningful[59].

---

[56] A positive weight function can be viewed as just a diagonal operator with real values on the diagonal, which is also therefore a Hermitian operator

[57] This name "transformed inner product" is one we are creating here.

[58] Note, incidentally, that, though transforms are often defined with unitary operators (see Eq. (140)) (or ones proportional to unitary operators, as in Fourier transforms, for example, there is no requirement that this operator $\mathsf{B}$ is unitary.

[59] In that case, we may want to know the SVD of the Green's function $\mathsf{G}_{SR}$ (which will be a Hilbert-Schmidt operator and hence compact) for the wave equation of interest when mapping from specific "source" Hilbert space to a specific "receiving" Hilbert space, for example. The resulting sets of functions will give us the "best" possible sources and corresponding received waves, all of which will be orthogonal in their respective spaces. These will also correspond

## 9.1    Form of singular value decomposition

For a compact (but not necessarily Hermitian) operator $\mathsf{A}$, mapping from a space $H_S$ to a possibly different space $H_R$, we consider first the eigenequation for the operator $\mathsf{A}^\dagger\mathsf{A}$, i.e.,

$$\mathsf{A}^\dagger\mathsf{A}\psi_j = c_j\psi_j \tag{126}$$

for the set of eigenvectors[60] $\{\psi_j\}$ (which we will choose to be normalized) in $H_S$ and the corresponding eigenvalues $c_j$. Then

$$c_j\left(\psi_j,\psi_j\right) = \left(\psi_j,\mathsf{A}^\dagger\mathsf{A}\psi_j\right) \equiv \left(\mathsf{A}\psi_j,\mathsf{A}\psi_j\right) \geq 0 \tag{127}$$

because in the last step in Eq. (127) we have an inner product of a vector with itself (see (IP4), Eq. (23)). So necessarily all the eigenvalues of $\mathsf{A}^\dagger\mathsf{A}$ (and similarly of $\mathsf{A}\mathsf{A}^\dagger$) are non-negative. So, we can choose to write these eigenvalues as $c_j = |s_j|^2$. So, using the expansion of the form Eq. (113) for $\mathsf{A}^\dagger\mathsf{A}$, we have

$$\mathsf{A}^\dagger\mathsf{A} = \sum_{j=1}^{\infty} |s_j|^2\, |\psi_j\rangle\langle\psi_j| \tag{128}$$

So,

$$\|\mathsf{A}\psi_n\|^2 = \langle\psi_n|\mathsf{A}^\dagger\mathsf{A}|\psi_n\rangle = |s_n|^2 \tag{129}$$

Then

$$\|\mathsf{A}\psi_n\| = |s_n| \tag{130}$$

So we can construct a set of functions $\{\phi_n\}$ in $H_R$ for all eigenfunctions $\{\psi_j\}$ corresponding to non-zero eigenvalues, where we define

$$|\phi_n\rangle = \frac{1}{s_n}\mathsf{A}|\psi_n\rangle \tag{131}$$

This set of functions is, first, normalized; that is

$$\langle\phi_n|\phi_n\rangle = \frac{1}{s_n^* s_n}\langle\psi_n|\mathsf{A}^\dagger\mathsf{A}|\psi_n\rangle = \frac{|s_n|^2}{s_n^* s_n} = 1 \tag{132}$$

and we have

$$
\begin{aligned}
\langle\phi_m|\phi_n\rangle &= \frac{1}{s_m^* s_n}\langle\psi_m|\mathsf{A}^\dagger\mathsf{A}|\psi_n\rangle = \frac{1}{s_m^* s_n}\langle\psi_m|\left(\sum_{j=1}^{\infty}|s_j|^2\,|\psi_j\rangle\langle\psi_j|\right)|\psi_n\rangle \\
&= \frac{1}{s_m^* s_n}\sum_{j=1}^{\infty}|s_j|^2\,\langle\psi_m|\psi_j\rangle\langle\psi_j|\psi_n\rangle = \frac{1}{s_m^* s_n}\sum_{j=1}^{\infty}|s_j|^2\,\langle\psi_m|\psi_j\rangle\delta_{jn} \\
&= \frac{|s_n|^2}{s_m^* s_n}\langle\psi_m|\psi_n\rangle = \frac{|s_n|^2}{s_m^* s_n}\delta_{mn} = \delta_{mn}
\end{aligned}
\tag{133}
$$

---

to the best-coupled and orthogonal channels for communicating with waves between the volumes [6]. The SVD approach also allows a way to synthesize arbitrary linear optical components [7].

[60] Note that these eigenvectors are orthogonal, being eigenvectors of a compact Hermitian operator, and with appropriately-chosen mutually orthogonal versions of any degenerate eigenvectors.

so this set $\{\phi_n\}$ is also orthonormal. Now suppose we consider an arbitrary function $\psi$ in $H_S$. Then we can expand it in the orthonormal set $\{\psi_j\}$ as in Eq. (81) to obtain

$$|\psi\rangle = \sum_j \langle \psi_j | \psi \rangle | \psi_j \rangle \tag{134}$$

So, using (131)

$$A|\psi\rangle = A\sum_j \langle \psi_j | \psi \rangle | \psi_j \rangle = \sum_j \langle \psi_j | \psi \rangle A | \psi_j \rangle = \sum_j s_j \langle \psi_j | \psi \rangle | \phi_j \rangle$$
$$= \sum_j s_j | \phi_j \rangle \langle \psi_j | \psi \rangle = \left( \sum_j s_j | \phi_j \rangle \langle \psi_j | \right) | \psi \rangle \tag{135}$$

Since $\psi$ was arbitrary, we can therefore write

$$A \equiv \sum_j s_j | \phi_j \rangle \langle \psi_j | \tag{136}$$

which is the *singular value decomposition* (often abbreviated to SVD) of the operator $A$ from a space $H_S$ to a possibly different space $H_R$. The numbers $s_j$ are called the *singular values* of the operator $A$.

Note, first, that we can perform this SVD for any compact operator. Second, this SVD tells us that we can view any such compact operator $A$ as "connecting" a set of orthogonal functions $\{\psi_j\}$ in $H_S$ one-by-one to a set of orthogonal functions $\{\phi_j\}$ in $H_R$, with associated "connection strengths" given by the corresponding singular values $s_j$ in each case.

## 9.2    Establishing the singular value decomposition of an operator

From Eq. (136), we can write

$$A^\dagger A = \sum_{j,k} \left( s_k^* | \psi_k \rangle \langle \phi_k | \right) \left( s_j | \phi_j \rangle \langle \psi_j | \right)$$
$$= \sum_{j,k} s_k^* s_j | \psi_k \rangle \langle \phi_k | \phi_j \rangle \langle \psi_j | = \sum_{j,k} \delta_{kj} s_k^* s_j | \psi_k \rangle \langle \psi_j | \tag{137}$$
$$= \sum_j |s_j|^2 | \psi_j \rangle \langle \psi_j |$$

But from Eq. (113), we see that this is just the representation of the operator on a basis of its eigenfunctions, which are $|\psi_j\rangle$ with eigenvalues $|s_j|^2$. Explicitly, we can check that these are the eigenfunctions and eigenvalues of $A^\dagger A$.

$$A^\dagger A | \psi_n \rangle = \left( \sum_{j,k} |s_j|^2 | \psi_j \rangle \langle \psi_j | \right) | \psi_n \rangle$$
$$= \sum_{j,k} |s_j|^2 | \psi_j \rangle \langle \psi_j | \psi_n \rangle = \sum_{j,k} |s_j|^2 | \psi_j \rangle \delta_{jn} \tag{138}$$
$$= |s_n|^2 | \psi_n \rangle$$

Similarly, the functions $|\phi_j\rangle$ are the eigenfunctions of $AA^\dagger$ with the same eigenvalues $|s_j|^2$. Hence the singular value decomposition can be established by solving for the eigenfunctions and eigenvalues of $A^\dagger A$ and for the eigenfunctions of $AA^\dagger$.

## 9.3    Matrix representation of the singular value decomposition

We can now explicitly show a common operator way of writing the SVD, which then also directly applies to matrices, and is then seen as a matrix factorization. Specifically, the SVD of some linear operator or matrix $A$ would be written in the form

$$A = VD_{diag}U^\dagger \tag{139}$$

where $U$ and $V$ are unitary operators or matrices, that is, operators or matrices for which

$$U^\dagger U = I_{op} \text{ and } V^\dagger V = I_{op} \tag{140}$$

where $I_{op}$ is the identity operator or matrix. We prove this equivalence below in 11.17 "Proof (17) of the equivalence of Dirac and matrix forms of SVD".

Another way of looking at this is that, if we expand $|\psi_p\rangle$ and $|\phi_q\rangle$ on some basis $\{\gamma_j\}$, then the elements of the $p$th row of $U^\dagger$ are the elements of $\langle\psi_p|$, and the elements of the $q$th column of $V$ are the elements of $|\phi_q\rangle$.

The SVD is, of course, a standard decomposition for finite matrices. Note here, though, that we are also rigorously defining the equivalent mathematics for compact operators that may be operating in or between infinite dimensional spaces.

# 10 Concluding remarks

This completes our introduction to this mathematics. Obviously, the reader can proceed further, and the various standard functional analysis texts certainly provide that. Indeed, my hope is that this introduction can make those texts[61] more accessible and hence valuable.

# 11 Proofs of specific mathematical theorems

## 11.1    Proof (1) that every convergent sequence in a metric space is a Cauchy sequence

If a sequence $(x_n)$ converges to the element $x$, i.e., $x_n \to x$, then for every (real number) $\varepsilon > 0$ there is a positive integer or natural number $N$ for any given $\varepsilon$ such that $d(x_n, x) < \varepsilon/2$ for all $n > N$. Hence by the triangle inequality for metrics and the symmetry of a metric, both of which are basic defined properties of a metric (see the definitions (M4) and (M3) respectively in the definitions (2) above),

$$d(x_m, x_n) \le d(x_m, x) + d(x, x_n) < (\varepsilon/2) + (\varepsilon/2) = \varepsilon$$

This shows that the sequence $(x_n)$ is Cauchy, completing the proof.

## 11.2    Proof (2) of the Bolzano-Weierstrass theorem

Here we prove the Bolzano-Weierstrass theorem (7) ("Each bounded sequence of real numbers has a convergent subsequence"). One standard proof uses "nested intervals". We can start by imagining that we have marked all the elements of an arbitrary (infinitely long) bounded sequence of real numbers $x_n$ (where

---

[61] Which mathematicians should understand are very difficult for ordinary mortals to follow!

*n* runs over all the natural numbers) on a real "line", as in Fig. 1(a); all the points $x_n$ necessarily lie between the lower bound, a number $x_{\text{inf}}$ corresponding to the infimum of the set of points, and the upper bound, a number $x_{\text{sup}}$ corresponding to the supremum of the set of points. Of course, the number of points we need to mark on the line is infinite, and for graphic purposes we can only indicate some of these on the graph, but we understand the actual number of points to be infinite.

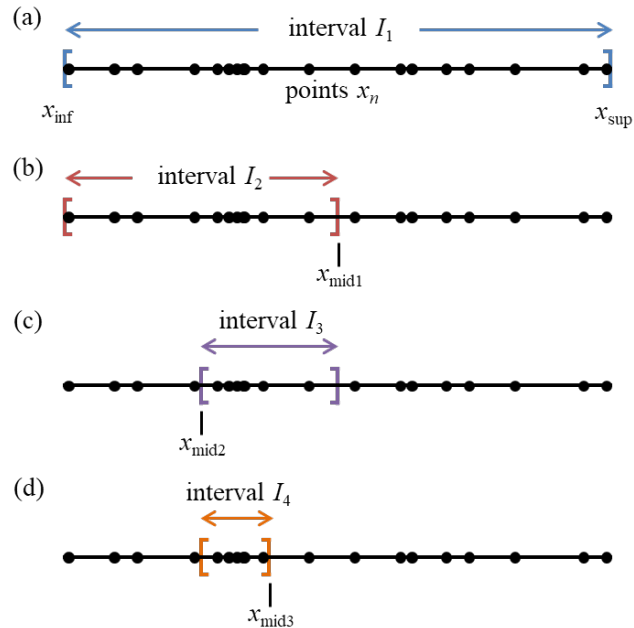

Fig. 1. Illustration of the process, starting with a sequence $(x_n)$ of points that are marked on the line, dividing an interval in two progressively, each time retaining an interval that has an infinite number of points, and hence contains an infinite subsequence of the original sequence $(x_n)$.

By definition, because we have an infinitely long sequence, then within the interval $I_1$, which goes from the infimum $x_{\text{inf}}$ to the supremum $x_{\text{sup}}$, there is an infinite number of points on the line. Now let us divide that interval in half, with a mid-point $x_{\text{mid1}}$. Our goal here is to establish a new interval, half as big as the previous one, and still with an infinite number of points in it. There are now three possibilities: (1) there is an infinite number of points in the interval between $x_{\text{inf}}$ and the mid-point $x_{\text{mid1}}$ but a finite number between $x_{\text{mid1}}$ and $x_{\text{sup}}$; (2) there is an infinite number of points in the interval between $x_{\text{mid1}}$ and $x_{\text{sup}}$ but a finite number between $x_{\text{inf}}$ the mid-point $x_{\text{mid1}}$; (3) there are infinite numbers of points between $x_{\text{inf}}$ and the mid-point $x_{\text{mid1}}$ as well as infinite number of points in the interval between $x_{\text{mid1}}$ and $x_{\text{sup}}$. In the first case, we now choose a new interval $I_2$ that runs from $x_{\text{inf}}$ and the mid-point $x_{\text{mid1}}$ (which is the example case shown in Fig. 1(b)). In the second case, we instead choose the new interval $I_2$ to run between $x_{\text{mid1}}$ and $x_{\text{sup}}$. In the third case, it does not matter which of the two new intervals we choose; we just arbitrarily choose one or the other; our goal is to show there is at least one convergent subsequence, so either one of these intervals would be suitable (it is not a problem if there are two convergent subsequences). The interval we are left with contains an (infinitely long) subsequence of the original sequence. (On whatever interval we end up choosing, we should choose it to include its end points so that we do not end up with a sequence that converges to a limit that lies "just" outside the interval).

Now we keep repeating this process, as illustrated in Fig. 1(c) and Fig. 1(d) for example successive intervals, dividing the interval in two each time, choosing a (or the) part with an infinite number of points within it, and continuing this process. As a result, we end up with an arbitrarily small interval that nonetheless contains a subsequence with an infinite number of points. Thus we can see we are establishing a convergent subsequence.

So, formally, after the choice of the $j$th interval, we have an (infinitely long) subsequence $y_{jm}$ (where $m$ runs over all the natural numbers) of the original sequence $x_n$. (Note that the elements of $y_{jm}$ are all elements of the original sequence $x_n$, and are in the same relative order as they were in $x_n$.) The size of this interval is $\Delta y_j = \left( x_{\sup} - x_{\inf} \right) / 2^{(j-1)}$ and all of its elements lie within this range (or on the edge of it). So, for any $\varepsilon$, no matter how small, there is always some sufficiently large choice of $j$ such that $\varepsilon < \Delta y_j$. Then, for our standard metric for real numbers $s$ and $t$, that is, $d\left( s, t \right) = |s - t|$, we have, for any elements $y_{jp}$ and $y_{jq}$, where $p$ and $q$ are any members of the set of natural numbers,

$$d\left( y_{jp}, y_{jq} \right) = \left| y_{jp} - y_{jq} \right| < \varepsilon \tag{141}$$

If we choose an $x$ that lies in the range between $y_{pj}$ and $y_{qj}$ (inclusive of the end points), then we can say for any $\varepsilon$, no matter how small, there is an $x$ such that

$$d\left( y_{jp}, x \right) = \left| y_{jp} - x \right| < \varepsilon \tag{142}$$

Hence, there is a convergent subsequence of the original sequence $\left( x_n \right)$ that approaches arbitrarily closely to some limit $x$, and so we formally have proved that for any bounded sequence $\left( x_n \right)$ there is a convergent subsequence, proving the theorem as required.

## 11.3   Proof (3) of the existence of a basis for a Hilbert space

If we can show that any vector $\gamma$ in a given Hilbert space $H$ can be expressed as a linear combination of orthonormal functions $\{ \alpha_1, \alpha_2, \ldots \}$, a set that may be of any required dimensionality (including infinite dimensionality), then by definition that set $\{ \alpha_1, \alpha_2, \ldots \}$ is a basis for the space $H$. To prove this set exists, we can formally construct it by considering as many non-zero vectors $\gamma_1, \gamma_2, \ldots$ as we like in the space, and showing how to construct the set of orthonormal functions $\{ \alpha_1, \alpha_2, \ldots \}$ from them. (This process is formally equivalent to *Gram-Schmidt orthogonalization*.)

We can start this process by taking the first of these vectors, $\gamma_1$, and constructing a normalized version of it, making that the first basis set element $\alpha_1$, i.e.,

$$\alpha_1 = \frac{\gamma_1}{\sqrt{\left( \gamma_1, \gamma_1 \right)}} \tag{143}$$

Now we consider the second vector $\gamma_2$. If $\gamma_2 = c\gamma_1$ where $c$ is some complex number, then we can already fully represent $\gamma_2$ already using just $\alpha_1$, i.e., explicitly $\gamma_2 = c\sqrt{\left( \gamma_1, \gamma_1 \right)} \alpha_1$. Since our goal is to establish the various elements $\alpha_1, \alpha_2, \ldots$ that will make up a basis set, then to save unnecessary operations, we will therefore just consider additional vectors $\gamma_2, \gamma_3, \ldots$ that cannot already be represented by the basis vectors we have already found. (Formally this means that we are choosing $\gamma_1, \gamma_2, \ldots$ to be linearly independent.). Therefore, with this restriction, we can write

$$\gamma_2 = \left( \gamma_2, \alpha_1 \right) \alpha_1 + \beta_2 \tag{144}$$

where $\beta_2$ is some non-zero vector orthogonal to $\alpha_1$. To see that $\beta_2$ is orthogonal, we can form the inner product

$$\begin{aligned} \left( \alpha_1, \gamma_2 \right) &= \left( \alpha_1, \gamma_2 \right)\left( \alpha_1, \alpha_1 \right) + \left( \alpha_1, \beta_2 \right) \\ &= \left( \alpha_1, \gamma_2 \right) + \left( \alpha_1, \beta_2 \right) \end{aligned} \tag{145}$$

so

$$\left( \alpha_1, \beta_2 \right) = 0 \tag{146}$$

proving the orthogonality. Now, therefore, we can form a second element of our basis set using a normalized version of $\beta_2$, specifically

$$\alpha_2 = \frac{\beta_2}{\sqrt{(\beta_2, \beta_2)}} \tag{147}.$$

To construct the third element, we choose a $\gamma_3$ that cannot already be represented as a linear combination of $\alpha_1$ and $\alpha_2$, leaving an orthogonal vector $\beta_3$ as in

$$\gamma_3 = \left(\sum_{j=1}^{2} (\alpha_j, \gamma_3) \alpha_j\right) + \beta_3 \tag{148}$$

Generally, we can keep going like this, with

$$\gamma_m = \left(\sum_{j=1}^{m-1} (\alpha_j, \gamma_m) \alpha_j\right) + \beta_m \tag{149}$$

and choosing

$$\alpha_m = \frac{\beta_m}{\sqrt{(\beta_m, \beta_m)}} \tag{150}$$

In this process, if our basis set is not complete, as proved by the fact that it cannot represent some vector, then we just add in a normalized version (i.e., $\alpha_m$) of the orthogonal "remainder" vector $\beta_m$ as the necessary new element in our basis set.

Of course, if we had a space of finite dimensionality, this process would truncate at some point once we could no longer find any vector in the space that could not be expressed as a linear combination of the basis vectors we had found so far, and we would have found our basis set. For an infinite dimensional space, we can just keep going, and so, inductively, we can create an orthonormal basis set to represent any function in such a Hilbert space.

## 11.4  Proof (4) of a criterion for compactness of an operator

We are proving the statement that "The operator $\mathsf{A}$ (from the normed space $F$ to the normed space $G$) is compact if and only if it maps every bounded sequence $(\alpha_m)$ of vectors in $F$ into a sequence in $G$ that has a convergent subsequence." This is mostly a question of definition[62].

We can deal with the "only if" part of this criterion first, which is true essentially directly from the definition (83) of a compact operator. Specifically, if the operator $\mathsf{A}$ is compact (as defined above in (83)) and the set of vectors $\alpha_m$ is bounded, then by that definition of compactness, the resulting set of vectors $\gamma_m = \mathsf{A}\alpha_m$ gives a space whose closure is compact.

Next we deal with the "if" part, which, restated, becomes "if $\mathsf{A}$ maps every bounded sequence $(\alpha_m)$ of vectors in $F$ into a sequence $(\gamma_m)$ (i.e., where $\gamma_m = \mathsf{A}\alpha_m$) in $G$ that has a convergent subsequence, then $\mathsf{A}$ is compact." Consider some subspace $B$ of $F$ that consists only of bounded vectors $\alpha_m$, and consider any sequence of vectors $(\gamma_m)$ in $G$ that can be formed by the action of $\mathsf{A}$ on any vector in $B$ (so the sequence $(\gamma_m)$ is an arbitrary sequence in the *image* of $B$ in $G$.) Then by presumption in the statement to be proved, this sequence $(\gamma_m)$ has a convergent subsequence. But $(\gamma_m)$ is an arbitrary sequence in the image of $B$, and $B$ is an arbitrary bounded subset of $F$, so every sequence in the image of $B$ has a convergent subsequence, so every sequence $(\gamma_m)$ generated by a bounded sequence $(\alpha_m)$ in $F$ has a bounded subsequence. So the resulting space in $G$ is compact, by the definition (8). Hence the operator $\mathsf{A}$ is compact by the definition (83). (The definition (83) is stated in terms of a precompact space, but a compact space is also precompact. A precompact space is one whose closure is compact. The closure of a compact space, which is already closed, is the same compact space.)

---

[62] This is one of those proofs that takes some space to write down, but that actually has very little in it.

## 11.5  Proof (5) of compactness of operators with finite dimensional range

Before proving this statement itself, we need to prove two other statements.

### 11.5.1  A theorem on linear combinations

First, we need to prove[63] the following theorem on linear combinations:

Consider a linearly independent set of vectors $\{\alpha_1, \ldots, \alpha_n\}$ in a normed vector space of any finite dimension. Then there is a real number $c > 0$ such that, for every choice of (possibly complex) scalars $a_1, \ldots, a_n$, we have

$$\|a_1 \alpha_1 + \cdots + a_n \alpha_n\| \geq c \left( |a_1| + \cdots |a_n| \right) \tag{151}$$

Before proving this theorem, we can note that, loosely, it is indicating that there are limits to how small a vector can be if it is made up out of linearly independent vectors that are large. The proof proceeds as follows.

We can write $s = |a_1| + \cdots |a_n|$ where $s \geq 0$ since the modulus of any complex number is greater than or equal to zero. The only way we can have $s = 0$ is for all the $a_j$ to be zero, in which case (151) holds for any $c$. So to complete the proof, we now consider all the other possibilities, for which necessarily $s > 0$. Then we can rewrite (151) as

$$\|b_1 \alpha_1 + \cdots + b_n \alpha_n\| \geq c \tag{152}$$

where $b_n = a_n / s$ and, necessarily, $\sum_{j=1}^{n} |b_j| = 1$. Hence, it is enough now to prove the existence of a $c > 0$ such that (152) holds for every collection of $n$ scalars $b_1, \ldots, b_n$ (complex numbers) with $\sum_{j=1}^{n} |b_j| = 1$.

Now we proceed by a *reductio ad absurdum* proof, starting by assuming that the statement is false, i.e., that there is a set or sets of such scalars for which $c$ is not greater than zero. To start this argument, we choose some (infinitely long) sequence $(\beta_m)$ of vectors, each of which can be written

$$\beta_m = b_1^{(m)} \alpha_1 + \cdots + b_n^{(m)} \alpha_n \tag{153}$$

with

$$\sum_{j=1}^{n} \left| b_j^{(m)} \right| = 1 \tag{154}$$

(the coefficients $b_j^{(m)}$ can be different for each such vector) and we require that this sequence is such that $\|\beta_m\| \to 0$ as $m \to \infty$.

Now we reason in what is sometimes called a "diagonal argument". Since $\sum_{j=1}^{n} \left| b_j^{(m)} \right| = 1$, we know that for every coefficient $b_j^{(m)}$ in any of the vectors $\beta_m$ in the sequence $(\beta_m)$

$$\left| b_j^{(m)} \right| \leq 1 \tag{155}$$

Since we have a sequence of vectors $(\beta_m)$, we can if we want construct a sequence of the values of the $j$th coefficient in each vector. Hence for each chosen $j$, we have a sequence of coefficients (a sequence of scalars, not of vectors)

$$\left( b_j^{(m)} \right) = \left( b_j^{(1)}, b_j^{(2)}, \ldots \right) \tag{156}$$

---

[63] This particular proof follows Kreyszig [2], Lemma 2.4-1.

If we imagined that we wrote out all the $n$ coefficients $b_1^{(1)},\ldots,b_n^{(1)}$ of the first vector $|\beta_1\rangle$ as a horizontal row, and then wrote the coefficients $b_1^{(2)},\ldots,b_n^{(2)}$ of the second vector $|\beta_2\rangle$ on a second horizontal row beneath it, and so on, as in

$$
\begin{array}{ccc}
b_1^{(1)} & \cdots & b_n^{(1)} \\
b_1^{(2)} & \cdots & b_n^{(2)} \\
\vdots & \vdots & \vdots
\end{array}
$$

then this sequence $\left(b_j^{(m)}\right) = \left(b_j^{(1)}, b_j^{(2)}, \ldots\right)$ would be one vertical column.

Note that, as usual with sequences, this sequence is infinitely long, and we know from (155) that this is a bounded sequence. Now we specifically choose the sequence $\left(b_1^{(m)}\right)$ (i.e., the first column). So, from the Bolzano-Weierstrass theorem, the sequence $\left(b_1^{(m)}\right)$ has a convergent subsequence, with some limit $b_1$.

Now we take the subsequence of vectors that corresponds to those with their first coefficient as this subsequence of $\left(b_1^{(m)}\right)$, with that first coefficient still limiting to $b_1$.

From that subsequence of vectors, we can choose a "sub-subsequence" (which is just another subsequence) in which the second coefficient similarly limits to some number $b_2$. (The existence and convergence of this (sub)subsequence is similarly guaranteed by the Bolzano-Weierstrass theorem.) We use this argument progressively a total of $n$ times, with each "column" converging to a corresponding limit $b_j$, by which time we are left with a (sub)sequence of vectors that we can call $(\gamma_k)$, a subsequence of the original sequence $(\beta_m)$. The individual vectors in the sequence $(\gamma_k)$ are of the form

$$\gamma_k = \sum_{j=1}^n g_j^{(k)} \alpha_j \tag{157}$$

and we have

$$\sum_{j=1}^n \left| g_j^{(k)} \right| = 1 \tag{158}$$

because the coefficients $g_j^{(k)}$, $j = 1,\ldots,n$ for each $k$ are just the coefficients $b_j^{(m)}$, $j = 1,\ldots,n$ for some vector $\beta_m$ in the original sequence of vectors; we have just been choosing a subsequence from that original sequence, and each $\gamma_k$ is just some $\beta_m$ in the original sequence. This sequence $(\gamma_k)$ converges to the vector

$$\gamma = \sum_{j=1}^n b_j \alpha_j \tag{159}$$

where, as in (154) and (158), we have

$$\sum_{j=1}^n \left| b_j \right| = 1 \tag{160}$$

so the $b_j$ cannot all be zero. Since the original vectors $\{\alpha_1,\ldots,\alpha_n\}$ were by choice linearly independent, then, with coefficients $b_j$ that are not all zero, the vector $\gamma$ cannot be the zero vector. We have found a convergent subsequence of the original sequence $(\beta_m)$ that does not converge to the zero vector. But this contradicts the original assumption that we could construct such a sequence of vectors that converges to zero, as required to allow the non-negative number $c$ to be zero. Hence, by *reductio ad absurdum, $c > 0$*, and we have proved the theorem.

### 11.5.2   A theorem on the compactness of any closed and bounded finite dimensional normed space

Next we need a second result, which is the following theorem[64]:

> In a finite-dimensional normed space $G$, any subset $M$ of that space is compact if it is closed and bounded.
> (161)

The proof of this theorem is as follows. We presume the space $G$ has dimensionality $n$ (so a linear combination of $n$ linearly independent vectors is required to specify an arbitrary vector in the space), and we presume that $\{\alpha_1, \ldots, \alpha_n\}$ is a basis for the space $G$. Now consider an infinitely long sequence of vectors $(\beta_m)$ in the subspace $M$. Each $\beta_m$ can be represented on the basis as

$$\beta_m = b_1^{(m)}\alpha_1 + \cdots + b_n^{(m)}\alpha_n \tag{162}$$

Since the subspace $M$ is bounded, so also is the sequence $(\beta_m)$, and we can call that bound some positive real number $v$, so for any $\beta_m$, we have $\|\beta_m\| \le v$. Then, using the result (151) just proved above, we have

$$v \ge \|\beta_m\| = \left\| \sum_{j=1}^{n} b_j^{(m)}\alpha_j \right\| \ge c \sum_{j=1}^{n} \left| b_j^{(m)} \right| \tag{163}$$

where $c > 0$. Hence the infinitely long sequence of numbers $(b_j^{(m)})$ for some fixed $j$ is bounded, and by the Bolzano-Weierstrass theorem, it must have an accumulation point $g_j$. By a similar "diagonal" argument as in the proof above of (151), we conclude that the infinitely long sequence $(\beta_m)$ has an infinitely long subsequence $(\gamma_k)$ that converges to a vector $\gamma = \sum_{j=1}^{n} b_j \alpha_j$. Since the (sub)space $M$ is closed, this vector $\gamma$ must be in the space $M$. Hence we have proved that the arbitrary (infinitely long) sequence $(\beta_m)$ has a subsequence that converges in $M$. Hence $M$ is compact, proving the theorem as in (161).

### 11.5.3   Core proof of theorem

Now that we have proved two key results, (151) and (161), we can proceed to the proof the theorem we want[65], which can be stated as follows:

> For a bounded operator $\mathsf{A}$ from a normed space $F$ to a normed space $G$, if the range of $\mathsf{A}$ in $G$ has finite dimensionality, then $\mathsf{A}$ is a compact operator.
> (164)

Note, incidentally, that it is only necessary that the range of the operator has finite dimension. Of course, a finite matrix has both a finite dimensional range and a finite dimensional domain, but the finite dimensional range is sufficient to prove that a (bounded) finite matrix is automatically a compact operator. We can prove the theorem (164) as follows.

Consider an arbitrary (infinitely long) bounded sequence $(\beta_m)$ of vectors in $F$. Then the inequality (49), which here becomes $\|\mathsf{A}\beta_m\| \le \|\mathsf{A}\|\|\beta_m\|$ , shows that the (infinitely long) sequence of vectors $(\mathsf{A}\beta_m)$ is bounded. Hence the set of vectors $\{\mathsf{A}\beta_m\}$ constitutes a bounded subset of $G$. Since this is a bounded set, by the Bolzano-Weierstrass theorem (on the sequence of real numbers formed using the metric for the vectors on the sequence), the sequence $(\mathsf{A}\beta_m)$ will have a convergent subsequence, converging to some vector $\gamma$. Hence the space generated by all the possible sequences $(\beta_m)$ is bounded, and can be closed by adding in the vectors corresponding to the accumulation points of these sequences, making the resulting space closed and bounded. But we are asserting that the space $G$ has finite dimensionality, and so this subspace must also have finite dimensionality. Since this subspace is a closed, bounded finite-

---

dimensionality space, it is therefore compact by (161). To close it, we had just to add the corresponding limiting vectors, and so the operator $\mathsf{A}$ is generating a precompact space when acting on bounded vectors, and it is therefore compact by the definition (83).

## 11.6 Proof (6) of compact limit operator for convergent compact operators

We remember that we are trying to prove the theorem (86)

> Consider an infinitely long sequence $(\mathsf{A}_n)$ of compact linear operators from a normed space $F$ into a Hilbert space $H$. If $\|\mathsf{A}_n - \mathsf{A}\| \to 0$ as $n \to \infty$ for some operator $\mathsf{A}$, then this limit operator $\mathsf{A}$ is compact.

This proof uses a "diagonal argument", which we introduced first above in the proof 11.5.1 "A theorem on linear combinations" (151) in 11.5 "Proof (5) of compactness of operators with finite dimensional range". In this way, we will show that for any bounded sequence $(\beta_m)$ in $F$, the "image" sequence $(\mathsf{A}\beta_m)$ in $H$ has a convergent subsequence, and hence by the condition (84) for compactness of an operator, the operator $\mathsf{A}$ is compact.

So, we proceed as follows. Since $\mathsf{A}_1$ (the first operator in the sequence $(\mathsf{A}_m)$) is compact, then it maps bounded sequences $(\beta_m)$ in $F$ to sequences $(\mathsf{A}_1\beta_m)$ that have a convergent subsequence in $H$. We notate that subsequence as $(\mathsf{A}_1\gamma_{1,m})$ for some corresponding sequence $(\gamma_{1,m})$ in $F$ that is a subsequence of $(\beta_m)$. Now a sequence that is convergent in a metric is also automatically a Cauchy sequence (see (6) above), a property we will use later, so instead of just saying that we have a convergent subsequence, we will say the subsequence is Cauchy. So, the subsequence $(\mathsf{A}_1\gamma_{1,m})$ is Cauchy. Now we can proceed in a "diagonal argument" fashion. Similarly, since the operator $\mathsf{A}_2$ is compact (and indeed all the operators $\mathsf{A}_n$ are compact by choice) we can find a subsequence of $(\gamma_{1,m})$, which we will call $(\gamma_{2,m})$ for which the sequence $(\mathsf{A}_2\gamma_{2,m})$ is Cauchy. Continuing in this fashion, we see that the "diagonal sequence" $(\eta_q) = (\gamma_{q,m})$ (where $q$ is a natural number) is a subsequence of $(\beta_m)$ such that, for every $n \le q$, $(\mathsf{A}_n\eta_q)$ is Cauchy. Now, by choice $(\beta_m)$ is bounded, and hence $(\eta_q)$ is bounded, say, $\|\eta_q\| \le c$ for some positive real $c$, for all $q$.

Having established these Cauchy sequences by this diagonal method, we can now proceed to use the presumed operator convergence ($\|\mathsf{A}_n - \mathsf{A}\| \to 0$ as $n \to \infty$) together with this Cauchy property.

Because $\|\mathsf{A}_n - \mathsf{A}\| \to 0$, there is an $n = p$ such that $\|\mathsf{A} - \mathsf{A}_p\| < \delta$ for any positive $\delta$ we choose. Specifically, we will choose to write $\delta = \varepsilon / 3c$ for some positive number $\varepsilon$. Since $(\mathsf{A}_n\eta_q)$ is Cauchy for every $q \ge n$, then there is a (natural number) $u \ge p$ such that

$$\|\mathsf{A}_p\eta_j - \mathsf{A}_p\eta_k\| < \frac{\varepsilon}{3} \text{ for all } j,k \ge u. \tag{165}$$

Now, suppose we have four vectors $\mu$, $\kappa$, $\rho$, and $\zeta$ in a normed vector space. Then we could write

$$\mu - \zeta = \mu - \kappa + \kappa - \rho + \rho - \zeta \tag{166}$$

So, by the triangle inequality for norms (property N4 in (1)), we could write

$$\begin{aligned}\|\mu - \zeta\| &\le \|\mu - \kappa\| + \|\kappa - \rho + \rho - \zeta\| \\ &\le \|\mu - \kappa\| + \|\kappa - \rho\| + \|\rho - \zeta\|\end{aligned} \tag{167}$$

So, similarly, we can write for $j,k \ge u$

$$\|A\eta_j - A\eta_k\| \le \|A\eta_j - A_p\eta_j\| + \|A_p\eta_j - A_p\eta_k\| + \|A_p\eta_k - A\eta_k\|$$

$$\le \|A - A_p\|\|\eta_j\| + \frac{\varepsilon}{3} + \|A_p - A\|\|\eta_k\| \tag{168}$$

$$< \frac{\varepsilon}{3c}c + \frac{\varepsilon}{3} + \frac{\varepsilon}{3c}c = \varepsilon$$

This shows that $\left(A\eta_q\right)$ is Cauchy and converges since $H$ is complete (being a Hilbert space). Hence, finally, for an arbitrary bounded sequence $\left(|\beta_m\rangle\right)$ in $F$, the sequence $\left(A\beta_m\right)$ has a convergent subsequence in $H$, and hence by the condition (84) for compactness of an operator, the operator $A$ is compact.

## 11.7  Proof (7) of equivalent statements of the Hilbert-Schmidt sum rule limit S

We start from the definition (87) that $S \equiv \sum_j \|A\alpha_j\|^2$ for a Hilbert-Schmidt operator $A$ operating on vectors in a Hilbert space $H_1$ with a complete basis $\{\alpha_1, \alpha_2, \ldots\}$ to generate vectors in a Hilbert space $H_2$. Note first that norm $\|A\alpha_j\|$ is a vector norm in Hilbert space $H_2$, and so $\|A\alpha_j\|^2$ can be written as an inner product in that space. Specifically

$$\|A\alpha_j\|^2 \equiv \left(\gamma_j, \gamma_j\right) \tag{169}$$

where $\gamma_j \equiv A\alpha_j$ is a vector in $H_2$. Now we can use the equivalence as in Eq. (42), so $\left(\gamma_j, \gamma_j\right) = \langle \gamma_j | \gamma_j \rangle$, the equivalence as in Eq. (63), and standard matrix-vector manipulations, giving

$$\left(\gamma_j, \gamma_j\right) = \left(\gamma_j, A\alpha_j\right) = \langle \gamma_j | A | \alpha_j \rangle = \left(|\gamma_j\rangle\right)^\dagger A | \alpha_j \rangle$$

$$= \left(A|\alpha_j\rangle\right)^\dagger A | \alpha_j \rangle = \langle \alpha_j | A^\dagger A | \alpha_j \rangle \tag{170}$$

Hence, the sum-rule limit can be rewritten as

$$S = \sum_j \langle \alpha_j | A^\dagger A | \alpha_j \rangle \equiv Tr(A^\dagger A) \tag{171}$$

where the notation on the right, $Tr(A^\dagger A)$, is a shorthand for the *trace* of the matrix, the trace being the sum of the diagonal elements of a matrix.

We can now prove three standard equivalences about Eq. (171), all of which are proved by introducing and/or eliminating versions of the identity operator or matrix for the space (as in Eq. (82)).

First, the trace of any matrix is independent of the (complete) basis used to represent it, so the result $S$ from Eq. (171) is the same no matter what the complete basis $\{|\alpha_j\rangle\}$ is. This is a standard result, but we give the proof here for completeness. We consider a second complete basis $\{|\beta_k\rangle\}$ on the space, so we have the identity operator, which we can write on this basis as $I_{op} = \sum_k |\beta_k\rangle\langle\beta_k|$ or on the $\{|\alpha_j\rangle\}$ basis as $I_{op} = \sum_j |\alpha_j\rangle\langle\alpha_j|$. So starting from the trace of an operator or matrix $B$ expressed on the $\{|\alpha_j\rangle\}$ basis, we proceed, introducing $I_{op}$ twice (with different summation indices), moving round complex numbers (inner products) and eliminating an identity operator, i.e.,

$$Tr(\mathsf{B}) = \sum_j \langle \alpha_j | \mathsf{B} | \alpha_j \rangle = \sum_j \langle \alpha_j | \mathsf{I}_{op} \mathsf{B} \mathsf{I}_{op} | \alpha_j \rangle = \sum_{j,k,p} \langle \alpha_j | \left( |\beta_k\rangle \langle \beta_k | \mathsf{B} | \beta_p \rangle \langle \beta_p | \right) | \alpha_j \rangle$$

$$= \sum_{j,k,p} \langle \alpha_j | \beta_k \rangle \langle \beta_k | \mathsf{B} | \beta_p \rangle \langle \beta_p | \alpha_j \rangle = \sum_{j,k,p} \langle \beta_k | \mathsf{B} | \beta_p \rangle \langle \beta_p | \alpha_j \rangle \langle \alpha_j | \beta_k \rangle$$

$$= \sum_{k,p} \langle \beta_k | \mathsf{B} | \beta_p \rangle \langle \beta_p | \left( \sum_j | \alpha_j \rangle \langle \alpha_j | \right) | \beta_k \rangle = \sum_{k,p} \langle \beta_k | \mathsf{B} | \beta_p \rangle \langle \beta_p | \mathsf{I}_{op} | \beta_k \rangle$$

$$= \sum_{k,p} \langle \beta_k | \mathsf{B} | \beta_p \rangle \langle \beta_p | \beta_k \rangle = \sum_{k,p} \langle \beta_k | \mathsf{B} | \beta_p \rangle \delta_{pk} = \sum_k \langle \beta_k | \mathsf{B} | \beta_k \rangle$$

(172)

Hence we have proved that the trace of an operator or matrix is independent of the basis used. Applying this to the result Eq. (171) allows us therefore to conclude that we get the same answer for $S$ independent of the (complete) basis used to evaluate it.

Second, introducing an identity operator $\mathsf{I}_{op} = \sum_j | \alpha_j \rangle \langle \alpha_j |$ inside the sum and using associativity of matrix-vector multiplication, we can show

$$S = \sum_j \langle \alpha_j | \mathsf{A}^\dagger \mathsf{I}_{op} \mathsf{A} | \alpha_j \rangle = \sum_j \langle \alpha_j | \mathsf{A}^\dagger \left( \sum_k | \alpha_k \rangle \langle \alpha_k | \right) \mathsf{A} | \alpha_j \rangle$$

$$= \sum_{j,k} \langle \alpha_j | \mathsf{A}^\dagger | \alpha_k \rangle \langle \alpha_k | \mathsf{A} | \alpha_j \rangle$$

$$= \sum_{j,k} a_{kj}^* a_{kj} = \sum_{j,k} |a_{kj}|^2$$

(173)

which is the sum of the modulus squared of all the matrix elements.

Third, starting from the middle line in Eq. (173),

$$S = \sum_{j,k} \langle \alpha_j | \mathsf{A}^\dagger | \alpha_k \rangle \langle \alpha_k | \mathsf{A} | \alpha_j \rangle$$

$$= \sum_{j,k} \langle \alpha_k | \mathsf{A} | \alpha_j \rangle \langle \alpha_j | \mathsf{A}^\dagger | \alpha_k \rangle = \sum_k \langle \alpha_k | \mathsf{A} \left( \sum_j | \alpha_j \rangle \langle \alpha_j | \right) \mathsf{A}^\dagger | \alpha_k \rangle$$

$$= \sum_k \langle \alpha_k | \mathsf{A} \mathsf{I}_{op} \mathsf{A}^\dagger | \alpha_k \rangle = \sum_k \langle \alpha_k | \mathsf{A} \mathsf{A}^\dagger | \alpha_k \rangle = Tr(\mathsf{A}\mathsf{A}^\dagger)$$

(174)

These three equivalences, Eqs. (172), (173), and (174) are the ones we set out to prove.

## 11.8  Proof (8) of the operator norm inequality for the Hilbert-Schmidt norm

We can write an arbitrary function $\eta$ in $H_1$ on this basis $\{\alpha_1, \alpha_2, \ldots\}$

$$| \eta \rangle = \sum_p h_p | \alpha_p \rangle$$

(175)

So, with a basis $\{\beta_1, \beta_2, \ldots\}$ in $H_2$ we can write for an operator $\mathsf{A}$ expanded as in Eq. (67)

$$\mathsf{A} | \eta \rangle = \sum_p h_p \mathsf{A} | \alpha_p \rangle = \sum_{j,k,p} a_{jk} h_p | \beta_j \rangle \langle \alpha_k | \alpha_p \rangle = \sum_{j,k} a_{jk} h_k | \beta_j \rangle$$

(176)

Now, by definition of the vector norm

$$\| \mathsf{A} \eta \|^2 = \left( \langle \eta | \mathsf{A}^\dagger \right) \left( \mathsf{A} | \eta \rangle \right) = \left( \sum_{p,q} \langle \beta_p | a_{pq}^* h_q^* \right) \left( \sum_{j,k} a_{jk} h_k | \beta_j \rangle \right) = \sum_{p,q,j,k} a_{pq}^* h_q^* a_{jk} h_k \delta_{pj}$$

$$= \sum_j \left[ \left( \sum_q a_{jq}^* h_q^* \right) \left( \sum_k a_{jk} h_k \right) \right] = \sum_j \left| \sum_k a_{jk} h_k \right|^2$$

$$\leq \sum_j \left( \sum_k |a_{jk}|^2 \sum_m |h_m|^2 \right) = \left( \sum_{j,k} |a_{jk}|^2 \right) \left( \sum_m |h_m|^2 \right) = \| \mathsf{A} \|_{HS}^2 \, \| \eta \|^2$$

(177)

(We have used the Cauchy-Schwarz inequality Eq. (208) in going from the second to third line – see 11.12 "Proof (12) of Cauchy-Schwarz inequality" below.) So, finally, we have proved, as required, that for the Hilbert-Schmidt norm,

$$\|A\eta\| \le \|A\|_{HS} \|\eta\| \tag{178}$$

## 11.9   Proof (9) of compactness of Hilbert-Schmidt operators

### 11.9.1   First step – an inequality for the Hilbert-Schmidt norm

First, we will establish that we can write an expression analogous to Eq. (49), but where the operator norm is the Hilbert-Schmidt norm rather than the supremum norm in that relation. That is, we want to prove, for a Hilbert-Schmidt operator $A$ and for any arbitrary vector $\mu$ in the Hilbert space $H$ on which it acts

$$\|A\mu\| \le \|A\|_{HS} \|\mu\| \tag{179}$$

First, let us define a vector $\beta_1$ that is a normalized version of the vector $\mu$, i.e.

$$\beta_1 = \frac{\mu}{\|\mu\|} \tag{180}$$

Quite generally, then, for the vector norm based on the inner product, it is straightforward that $\|A\mu\| = \|A\beta_1\| \|\mu\|$. So now to prove Eq. (179), we need to prove that

$$\|A\beta_1\| \le \|A\|_{HS} \tag{181}$$

or equivalently

$$\|A\beta_1\|^2 \equiv \langle \beta_1 | A^\dagger A | \beta_1 \rangle \le \|A\|_{HS}^2 \tag{182}$$

Now, we are free to choose $\beta_1$ to be the first element of an orthogonal set that forms a basis for this Hilbert space $H$ of interest. So we have from Eq. (89)

$$\|A\|_{HS}^2 = \sum_k \langle \beta_k | A^\dagger A | \beta_k \rangle \ge \langle \beta_1 | A^\dagger A | \beta_1 \rangle \tag{183}$$

because all the elements $\langle \beta_k | A^\dagger A | \beta_k \rangle$ in the sum over $k$ are greater than or equal to zero, being inner products of the vector $A | \beta_k \rangle$ with itself. Hence

$$\|A\|_{HS}^2 \ge \langle \beta_1 | A^\dagger A | \beta_1 \rangle \equiv \|A\beta_1\|^2 \tag{184}$$

proving (181), and hence proving Eq. (179).

### 11.9.2   Second step – the remaining proof

For the next step in the larger proof, first we write a Hilbert-Schmidt operator $A$, which maps from a Hilbert space $H_1$ to a Hilbert space $H_2$, for arbitrary orthonormal basis sets $\left\{ |\alpha_k\rangle \right\}$ and $\left\{ |\beta_j\rangle \right\}$ in $H_1$ and $H_2$ respectively, in the form

$$A \equiv \sum_{j,k} a_{jk} | \beta_j \rangle \langle \alpha_k | \tag{185}$$

as in Eq. (65) (which we can always do). Next, we now consider another operator $A_n$ of in which we truncate the sum over one of the indices so that its range has finite dimensionality, that is,

$$A_n = \sum_{j=1}^{n} \sum_k a_{jk} | \beta_j \rangle \langle \alpha_k | \tag{186}$$

We note immediately that such an operator is compact, as proved above in 11.5 "Proof (5) of compactness of operators with finite dimensional range". Now consider the operator $A - A_n$, which, from Eqs. (185) and (186), we can write as

$$A - A_n = \sum_{j=n+1}^{\infty} \sum_k a_{jk} |\beta_j\rangle \langle \alpha_k|$$ (187)

So, from Eq. (89),

$$\|A - A_n\|_{HS}^2 = \sum_{j=n+1}^{\infty} \sum_k |a_{jk}|^2$$ (188)

But because we know $\sum_{j,k} |a_{jk}|^2$ is bounded because $A$ is a Hilbert-Schmidt operator, then

$$\|A - A_n\| \to 0 \text{ as } n \to \infty$$ (189)

Hence from the theorem (86), since $A$ is then the limit of a sequence of compact operators, $A$ is also compact. Hence we have proved our result that all Hilbert-Schmidt operators are compact.

## 11.10 Proof (10) of approximation of Hilbert-Schmidt operators by sufficiently large matrices

This proof follows similar approaches to the 11.9 "Proof (9) of compactness of Hilbert-Schmidt operators" above. For a matrix representation of the Hilbert-Schmidt operator $A$ that maps from Hilbert space $H_1$ to Hilbert space $H_2$ (which can be but is not necessarily different from $H_1$), written as in Eq. (185), we can write the "truncated" $m \times n$ matrix version $A_{mn}$ as

$$A_{mn} = \sum_{j=1}^{m} \sum_{k=1}^{n} a_{jk} |\beta_j\rangle \langle \alpha_k|$$ (190)

For some arbitrary (finite) vector $\mu$ in $H_1$, consider the vector $\eta$ (in $H_2$) that is the difference between the vectors $A\mu$ and $A_{mn}\mu$, i.e.,

$$\eta = A\mu - A_{mn}\mu = (A - A_{mn})\mu$$ (191)

Then

$$\|\eta\| = \|(A - A_{mn})\mu\| \le \|A - A_{mn}\|_{HS} \|\mu\|$$ (192)

where we have used the result Eq. (179) proved above in 11.9 "Proof (9) of compactness of Hilbert-Schmidt operators". So, from Eq. (89),

$$\|A - A_{mn}\|_{HS}^2 = \sum_{j=m+1}^{\infty} \sum_{k=n+1}^{\infty} |a_{jk}|^2$$ (193)

But because we know $\sum_{j,k} |a_{jk}|^2$ is bounded because $A$ is a Hilbert-Schmidt operator, then

$$\|A - A_{mn}\| \to 0 \text{ as } n \to \infty$$ (194)

and so $\|\eta\| \to 0$ as $m$ and $n$ tend to infinity. Because this difference vector vanishes in this limit, we have proved that we can approximate any Hilbert-Schmidt operator arbitrarily well by a sufficiently large matrix.

## 11.11 Proof (11) of Hilbert-Schmidt and compact nature of various operators derived from Hilbert-Schmidt operators

We now prove the Hilbert-Schmidt and compact nature of the operators $A^\dagger$, $A^\dagger A$ and $AA^\dagger$ for a Hilbert-Schmidt operator $A$. Consider a Hilbert-Schmidt operator $A$ operating on functions or vectors in a Hilbert space $H_1$ to generate functions or vectors in a Hilbert space $H_2$. Then from Eq. (89) we know $\sum_{j,k}|a_{kj}|^2$ is bounded (indeed, this boundedness is a necessary and sufficient condition for the corresponding operator to be Hilbert-Schmidt). For such matrix elements $a_{kj}$ of the operator $\mathsf{A}$, then we know from Eq. (76) that the matrix elements of the operator $\mathsf{A}^\dagger$ are $b_{jk} = a_{kj}^*$. Hence $\sum_{j,k}|b_{jk}|^2 = \sum_{j,k}|a_{kj}|^2$, which is therefore also bounded, and so $\mathsf{A}^\dagger$ is also a Hilbert-Schmidt operator (and therefore is also compact).

Given that $\sum_{j,k}|a_{kj}|^2$ is bounded, then for any matrix element $|a_{kj}|$ is also bounded, so for some sufficiently large positive real number $c$

$$|a_{kj}| \leq c \tag{195}$$

Now, by definition, and using Eq. (186) to represent $\mathsf{A}$, we can similarly write

$$\mathsf{A}^\dagger = \sum_{j,k} a_{jk}^* |\alpha_k\rangle \langle \beta_j| \tag{196}$$

and hence we can write the matrix product

$$
\begin{aligned}
\mathsf{A}^\dagger \mathsf{A} &= \sum_{j,k,p,q} a_{pq}^* |\alpha_q\rangle \langle \beta_p| a_{jk} |\beta_j\rangle \langle \alpha_k| \\
&= \sum_{j,k,p,q} a_{pq}^* |\alpha_q\rangle \delta_{pj} a_{jk} \langle \alpha_k| \\
&= \sum_{k,q} \left( \sum_j a_{jq}^* a_{jk} \right) |\alpha_q\rangle \langle \alpha_k|
\end{aligned}
\tag{197}
$$

So,

$$\|\mathsf{A}^\dagger \mathsf{A}\|_{HS}^2 \equiv \sum_{k,q} \left| \sum_j a_{jq}^* a_{jk} \right|^2 \tag{198}$$

Now

$$\left| \sum_j a_{jq}^* a_{jk} \right| \leq \sum_j \left| a_{jq}^* a_{jk} \right| \leq \sqrt{\sum_j |a_{jq}|^2} \sqrt{\sum_p |a_{pk}|^2} \tag{199}$$

where we used the Cauchy-Schwarz inequality for the last step. (We prove this inequality below in 11.12 "Proof (12) of Cauchy-Schwarz inequality".) So

$$
\begin{aligned}
\|\mathsf{A}^\dagger \mathsf{A}\|_{HS}^2 &\equiv \sum_{k,q} \left| \sum_j a_{jq}^* a_{jk} \right|^2 \leq \sum_{k,q} \left( \sum_j |a_{jq}|^2 \sum_p |a_{pk}|^2 \right) = \sum_{k,p} \left( \sum_{q,j} |a_{jq}|^2 \right) |a_{pk}|^2 \\
&= \sum_{k,p} \|\mathsf{A}\|_{HS}^2 |a_{pk}|^2 = \|\mathsf{A}\|_{HS}^2 \sum_{k,p} |a_{pk}|^2 = \|\mathsf{A}\|_{HS}^2 \|\mathsf{A}\|_{HS}^2 = \|\mathsf{A}\|_{HS}^4
\end{aligned}
\tag{200}
$$

So

$$\|\mathsf{A}^\dagger \mathsf{A}\|_{HS} \leq \|\mathsf{A}\|_{HS}^2 \tag{201}$$

Because $\mathsf{A}$ by choice is a Hilbert-Schmidt operator, it has a finite Hilbert-Schmidt norm $\|\mathsf{A}\|_{HS}$, and so $\|\mathsf{A}\|_{HS}^2$ is also finite. Hence $\|\mathsf{A}^\dagger \mathsf{A}\|_{HS}$ is finite, and so $\mathsf{A}^\dagger \mathsf{A}$ is a Hilbert-Schmidt operator, and hence also is compact. Since the operator $\mathsf{A}\mathsf{A}^\dagger$ is just the Hermitian adjoint of the operator $\mathsf{A}^\dagger \mathsf{A}$, it also is a Hilbert-Schmidt operator and is also compact.

## 11.12 Proof (12) of Cauchy-Schwarz inequality

Consider two mathematical vectors $|\beta\rangle$ and $|\gamma\rangle$ with complex elements. We can write these as column vectors if we wish

$$|\beta\rangle = \begin{bmatrix} b_1 \\ b_2 \\ \vdots \end{bmatrix} \text{ and } |\gamma\rangle = \begin{bmatrix} g_1 \\ g_2 \\ \vdots \end{bmatrix} \tag{202}$$

We presume these vectors are non-zero (the resulting theorem is trivial of either of them is zero). We also presume these two vectors are not proportional to one another (i.e., they are not in the same "direction") - again, the resulting theorem is trivially obvious if they are. Now we define a third vector

$$|\eta\rangle = |\beta\rangle - \frac{\langle\gamma|\beta\rangle}{\langle\gamma|\gamma\rangle}|\gamma\rangle \tag{203}$$

Here, the notation $\langle\gamma|\beta\rangle$ is just signifying that we are taking the simple Cartesian inner product of these vectors. For $|\eta\rangle$ to be zero, we would require $|\beta\rangle \propto |\gamma\rangle$, which by assumption is not the case, so $|\eta\rangle$ is non-zero.

Now we can write

$$\langle\gamma|\eta\rangle = \langle\gamma|\beta\rangle - \frac{\langle\gamma|\beta\rangle}{\langle\gamma|\gamma\rangle}\langle\gamma|\gamma\rangle = \langle\gamma|\beta\rangle - \langle\gamma|\beta\rangle = 0 \tag{204}$$

where we have used the inner product properties (IP1) and (IP2). (Since both $|\gamma\rangle$ and $|\eta\rangle$ are non-zero, for such an inner product to be zero, $|\eta\rangle$ is necessarily orthogonal to $|\gamma\rangle$.)

Rewriting Eq. (203) gives

$$|\beta\rangle = |\eta\rangle + \frac{\langle\gamma|\beta\rangle}{\langle\gamma|\gamma\rangle}|\gamma\rangle \tag{205}$$

so

$$\begin{aligned}
\|\beta\|^2 &\equiv \langle\beta|\beta\rangle = \left|\frac{\langle\gamma|\beta\rangle}{\langle\gamma|\gamma\rangle}\right|^2 \langle\gamma|\gamma\rangle + \langle\eta|\eta\rangle + \frac{\langle\gamma|\beta\rangle^*}{\langle\gamma|\gamma\rangle}\langle\gamma|\eta\rangle + \frac{\langle\gamma|\beta\rangle}{\langle\gamma|\gamma\rangle}\langle\eta|\gamma\rangle \\
&= \left|\frac{\langle\gamma|\beta\rangle}{\langle\gamma|\gamma\rangle}\right|^2 \langle\gamma|\gamma\rangle + \langle\eta|\eta\rangle \geq \left|\frac{\langle\gamma|\beta\rangle}{\langle\gamma|\gamma\rangle}\right|^2 \langle\gamma|\gamma\rangle = \frac{|\langle\gamma|\beta\rangle|^2}{\langle\gamma|\gamma\rangle} = \frac{|\langle\gamma|\beta\rangle|^2}{\|\gamma\|^2}
\end{aligned} \tag{206}$$

where we used Eq. (204) to eliminate the two right-most terms on the top line. So

$$|\langle\gamma|\beta\rangle| \leq \|\beta\|\|\gamma\| \tag{207}$$

or, equivalently, in component or summation form

$$\left|\sum_j g_j^* b_j\right| \leq \sqrt{\sum_p |g_p|^2}\sqrt{\sum_q |b_q|^2} \tag{208}$$

Eqs. (207) and (208) are each forms of the *Cauchy-Schwarz inequality*, which we had set out to prove.

## 11.13 Proof (13) of finite multiplicity

Consider a compact Hermitian operator $\mathsf{A}$ that acts on vectors in a Hilbert space $H$. If an eigenvalue $c$ of this operator has a multiplicity (degeneracy) different from 1, then there is a set of linearly independent vectors $\{\eta_j\}$ in $H$ of dimensionality given by the multiplicity. From those vectors we can form an orthonormal basis set $\{\beta_j\}$ of that same multiplicity that spans the space including all the vectors $\{\eta_j\}$.

That basis set is also a set of eigenvectors of $\mathsf{A}$ with the same eigenvalue $c$ since each element of the set is a linear combination of vectors that are eigenvectors with the same eigenvalue $c$. So, the effect of $\mathsf{A}$ on each orthogonal vector $\beta_j$ is to generate a vector, $c\beta_j$. Suppose $c$ is non-zero, and suppose that its multiplicity is infinite, so $\{\beta_j\}$ is an infinite set that is mapped, vector by vector, to the set $\{c\beta_j\}$, all of which vectors have the same finite non-zero norm and all of which are orthogonal. So, $\mathsf{A}$ maps the infinite sequence $(\beta_j)$ to the infinite sequence $(c\beta_j)$. But, because the vectors $\beta_j$ are orthonormal by definition, the sequence $(c\beta_j)$ has no convergent subsequence; the metric[66] $d(c\beta_n, c\beta_m) = \sqrt{2}|c|$ for any choice of two different values of $n$ and $m$. This contradicts the requirement for a compact operator ((84)) that it should map any sequence of finite vectors to a sequence with a convergent subsequence. Hence:

> For any non-zero eigenvalue $c$ of a compact Hermitian operator, the multiplicity of the eigenvalue is finite.

(209)

## 11.14 Proof (14) that the eigenvalues of Hermitian operator on an infinite dimensional space tend to zero

Suppose that we have a compact infinite-dimensional Hermitian operator that acts on vectors in a Hilbert space $H$. Consider the corresponding infinite set of eigenvectors $\{\alpha_j\}$, all of which are orthogonal (or can always be chosen to be so in the case of eigenvectors of degenerate eigenvalues.) Consider now the sequence $(\alpha_j)$ formed from those eigenvectors, in which each eigenvector appears only once, and in which we can consider the eigenvectors to be ordered in decreasing order of their corresponding eigenvalue. (The order within the finite set of eigenvectors corresponding to any degenerate eigenvalue does not matter.) If the eigenvalues do not converge to zero with increasing $j$, the resulting sequence of vectors $(\mathsf{A}\alpha_j)$ will not have a convergent subsequence, because they are all orthogonal vectors of finite length (see the "extreme example" argument above in "An illustrative extreme example" in the discussion of compact operators). Hence:

> If a compact Hermitian operator is operating on an infinite dimensional space, then the sequence of eigenvalues $(c_p)$ must tend to zero as $p \to \infty$.

(210)

## 11.15 Proof (15) of Hermitian operator supremum norm

We want to prove Eq. (110), $\|\mathsf{A}\| = \sup_{\|\alpha\|=1}\left|(\alpha, \mathsf{A}\alpha)\right|$, starting from Eq. (109), $\|\mathsf{A}\| = \sup_{\|\alpha\|=1}\|\mathsf{A}\alpha\|$ where $\alpha$ is any vector in the Hilbert space $H$ on which the Hermitian operator $\mathsf{A}$ acts. We do this by proving, first, that $\sup_{\|\alpha\|=1}\left|(\alpha, \mathsf{A}\alpha)\right| \leq \|\mathsf{A}\|$ and, second, that $\|\mathsf{A}\| \leq \sup_{\|\alpha\|=1}\left|(\alpha, \mathsf{A}\alpha)\right|$, which leaves the equality as the only remaining possibility from the two inequalities, which is what we want to prove.

For algebraic convenience, we define the (positive) real number $u$ as

$$u = \sup_{\|\alpha\|=1}\left|(\alpha, \mathsf{A}\alpha)\right|$$

(211)

(Note that this is equivalent to a statement $u = \sup_{\|\beta\|=1}\left|(\beta, \mathsf{A}\beta)\right|$; the name of the vector being used to find the supremum is arbitrary, and we will need this flexibility below.)

Using the Cauchy-Schwarz inequality, as in Eq. (207), and noting that $\|\alpha\| = 1$ by choice, we note, then, that

$$\left|(\alpha, \mathsf{A}\alpha)\right| \leq \|\mathsf{A}\alpha\|\|\alpha\| = \|\mathsf{A}\alpha\|$$

(212)

---

[66] Remember that the distance between the tips of two orthogonal unit vectors is $\sqrt{2}$.

So,

$$u = \sup_{\|\alpha\|=1} \left| (\alpha, A\alpha) \right| \leq \sup_{\|\alpha\|=1} \|A\alpha\| = \|A\| \tag{213}$$

(where the last step on the right is just the definition of $\|A\|$) completing the first half of the proof.

For the second half of the proof, note first that we can prove, for any vector $\eta$ in $H$, and considering all vectors $\gamma$ in $H$ with unit norm, i.e., $\|\gamma\| = 1$,

$$\|\eta\| = \sup_{\|\gamma\|=1} \left| (\gamma, \eta) \right| \tag{214}$$

To prove this statement, Eq. (214), note first that, by the Cauchy-Schwarz inequality, Eq. (207)

$$\left| (\gamma, \eta) \right| \leq \|\gamma\| \|\eta\| = \|\eta\| \tag{215}$$

so $\|\eta\| \geq \left| (\gamma, \eta) \right|$. But if we choose $\gamma = \eta / \|\eta\|$ (i.e., a normalized version of $\eta$), then

$$(\gamma, \eta) = (\eta, \eta) / \|\eta\| = \|\eta\| \tag{216}$$

since $\|\eta\| = \sqrt{(\eta, \eta)}$ by definition, which shows there is at least one choice of $\gamma$ for which $\|\eta\| = \left| (\gamma, \eta) \right|$. Hence, taking this result together with $\|\eta\| \geq \left| (\gamma, \eta) \right|$ from (215) proves Eq. (214).

So, with the definition Eq. (109), $\|A\| = \sup_{\|\alpha\|=1} \|A\alpha\|$, and choosing $\eta = A\alpha$ in Eq. (214)

$$\|A\| = \sup_{\|\alpha\|=1} \|A\alpha\| = \sup_{\|\alpha\|=1} \left( \sup_{\|\gamma\|=1} \left| (\gamma, A\alpha) \right| \right) \equiv \sup_{\|\alpha\|=1, \|\gamma\|=1} \left| (\gamma, A\alpha) \right| \tag{217}$$

Next we need to derive an inequality for $\left| (\gamma, A\alpha) \right|$. The first step is to prove an algebraic equivalence, and we start by choosing a vector $\mu = \exp(is)\gamma$ where the real number $s$ is chosen so that $(\mu, A\alpha)$ is real. (We are always free to do this, and such a number $s$ can always be found.) We note that

$$\left| (\mu, A\alpha) \right| = \left| (\exp(-is)\gamma, A\alpha) \right| = \left| \exp(is)(\gamma, A\alpha) \right| = \left| (\gamma, A\alpha) \right| \tag{218}$$

and

$$\|\mu\| = \sqrt{\langle \mu | \mu \rangle} = \sqrt{\langle \exp(-is)\gamma | \exp(-is)\gamma \rangle} = \sqrt{\langle \gamma | \gamma \rangle} = \|\gamma\| \tag{219}$$

Then we note that

$$
\begin{aligned}
&\left( \alpha + \mu, A(\alpha + \mu) \right) + \left( \alpha - \mu, A(\alpha - \mu) \right) \\
&= (\alpha, A\alpha) + (\mu, A\alpha) + (\alpha, A\mu) + (\mu, A\mu) \\
&\quad - (\alpha, A\alpha) + (\mu, A\alpha) + (\alpha, A\mu) - (\mu, A\mu) \\
&= 2\left[ (\mu, A\alpha) + (\alpha, A\mu) \right] \\
&= 2\left[ (\mu, A\alpha) + (A^\dagger \alpha, \mu) \right] \text{ by the definition of } A^\dagger \\
&= 2\left[ (\mu, A\alpha) + (A\alpha, \mu) \right] \text{ by the Hermiticity of A} \\
&= 2\left[ (\mu, A\alpha) + (\mu, A\alpha)^* \right] \text{ by (IP3)} \\
&= 4\text{Re}(\mu, A\alpha) \\
&= 4(\mu, A\alpha) \text{ by the chosen reality of this inner product}
\end{aligned}
\tag{220}
$$

So, from Eqs. (218) and (220)

$$\left| (\gamma, A\alpha) \right|^2 = \frac{1}{16} \left| \left( \alpha + \mu, A(\alpha + \mu) \right) + \left( \alpha - \mu, A(\alpha - \mu) \right) \right|^2 \tag{221}$$

Writing $\phi = \alpha + \mu$ and $\psi = \alpha - \mu$, we have

$$\left| (\gamma, \mathsf{A}\alpha) \right|^2 = \frac{1}{16} \left| (\phi, \mathsf{A}\phi) + (\psi, \mathsf{A}\psi) \right|^2 \tag{222}$$

Now

$$\left| (\phi, \mathsf{A}\phi) + (\psi, \mathsf{A}\psi) \right| \le \left| (\phi, \mathsf{A}\phi) \right| + \left| (\psi, \mathsf{A}\psi) \right| \tag{223}$$

by the triangle inequality for complex numbers. Writing a normalized vector $\beta = \phi / \|\phi\|$, and using the Cauchy-Schwarz inequality and the definition of $u$, Eq. (211),

$$\left| (\phi, \mathsf{A}\phi) \right| = \|\phi\|^2 \left( \beta, \mathsf{A}\beta \right) \le \|\phi\|^2 \sup_{\|\beta\|=1} \left| (\beta, \mathsf{A}\beta) \right| = \|\phi\|^2 \, u \tag{224}$$

and similarly

$$\left| (\psi, \mathsf{A}\psi) \right| \le \|\psi\|^2 \, u \tag{225}$$

So, using these results (224) and (225) in (223) and substituting that result into Eq. (222), we have

$$
\begin{aligned}
\left| (\gamma, \mathsf{A}\alpha) \right|^2 &\le \frac{u^2}{16} \left| \|\phi\|^2 + \|\psi\|^2 \right|^2 \\
&= \frac{u^2}{16} \left| \|\alpha + \mu\|^2 + \|\alpha - \mu\|^2 \right|^2 \equiv \frac{u^2}{16} \left| (\alpha + \mu, \alpha + \mu) + (\alpha - \mu, \alpha - \mu) \right|^2 \\
&= \frac{u^2}{16} \left| 2(\alpha, \alpha) + 2(\mu, \mu) \right|^2 \equiv \frac{u^2}{16} \left| 2\|\alpha\|^2 + 2\|\mu\|^2 \right|^2 \\
&\equiv \frac{u^2}{4} \left| \|\alpha\|^2 + \|\mu\|^2 \right|^2 = \frac{u^2}{4} \left| \|\alpha\|^2 + \|\gamma\|^2 \right|^2 = \frac{u^2}{4} |1 + 1|^2 \\
&= u^2
\end{aligned}
\tag{226}
$$

In the last step we used the fact that both $\|\alpha\|$ and $\|\gamma\|$ are 1, by choice in this proof. (The equivalence $\|\alpha + \mu\|^2 + \|\alpha - \mu\|^2 = 2\|\alpha\|^2 + 2\|\mu\|^2$ that is proved in the middle of these steps is sometimes called the *parallelogram law* or *parallelogram inequality*.) Hence, we have proved

$$\left| (\gamma, \mathsf{A}\alpha) \right| \le u \tag{227}$$

and so, automatically, $\sup_{\|\alpha\|=1, \|\gamma\|=1} \left| (\gamma, \mathsf{A}\alpha) \right| \le u$, and hence, using Eq. (217)

$$\|\mathsf{A}\| \le u \tag{228}$$

Since we have proved both that $\|\mathsf{A}\| \le u$ (Eq. (228)), and that $\|\mathsf{A}\| \ge u$ (Eq. (213), then we conclude that $\|\mathsf{A}\| = u = \sup_{\|\alpha\|=1} \left| (\alpha, \mathsf{A}\alpha) \right|$, which is the statement, Eq. (110), that we set out to prove.

## 11.16 Proof (16) of the spectral theorem

To prove the spectral theorem[67], first we prove for our compact, Hermitian operator $\mathsf{A}$ that either $r = \|\mathsf{A}\|$ or $r = -\|\mathsf{A}\|$ is an eigenvalue of $\mathsf{A}$.

---

[67] Our proof here is similar to that in [3], theorem 9.16, pp. 225 – 227, though our version is expanded. The overall structure of this proof is standard, and similar proofs are found in many other sources.

### 11.16.1 Proof that either $r = \|A\|$ or $r = -\|A\|$ is an eigenvalue of A

If $A = 0$ (the zero operator), then there is nothing to prove because all vectors in $H$ are then eigenvectors with eigenvalue zero, so we presume $A \neq 0$. Starting from the result $\|A\| = \sup_{\|\alpha\|=1} \left| (\alpha, A\alpha) \right|$ (Eq. (110), proved above in 11.15 "Proof (15) of Hermitian operator supremum norm"), we conclude that there must therefore be a sequence $(\alpha_n)$ in $H$, with $\|\alpha_n\| = 1$, such that

$$\|A\| = \lim_{n \to \infty} \left| (\alpha_n, A\alpha_n) \right| \tag{229}$$

Since $(\alpha_n, A\alpha_n)$ is necessarily real, as proved above (Eq. (99), (OE1)), then

$$\lim_{n \to \infty} (\alpha_n, A\alpha_n) = r \tag{230}$$

where, therefore, the real number

$$r = \|A\| \text{ or } r = -\|A\| \tag{231}$$

The (infinitely long) sequence $(\alpha_n)$ is bounded, and $A$ is compact. So there is a (infinitely long) subsequence $(\alpha_m)$, such that the sequence $(A\alpha_m)$ converges to some limit vector $\gamma$. That is,

$$\lim_{m \to \infty} A\alpha_m \to \gamma \tag{232}$$

Note that for this infinitely long sequence $(\alpha_m)$, because it is a subsequence of $(\alpha_n)$, we must still have, as in Eq. (230)

$$\lim_{m \to \infty} (\alpha_m, A\alpha_m) = r \tag{233}$$

Next, we will prove that $\gamma$ is an eigenvector of $A$, with eigenvalue $r$. First, we note that $\gamma$ is not the zero vector because then, from Eq. (230), we would have $r = 0$, and that cannot be because $r = \|A\|$ or $r = -\|A\|$ and by presumption $A \neq 0$. Now, if and only if $\gamma$ is an eigenvector with eigenvalue $r$, then $A\gamma = r\gamma$ (Eq. (100), (OE1)) so, formally, with the identity operator $I_{op}$ for the space $H$,

$$\left( A - r I_{op} \right) \gamma = 0 \text{ (the zero vector)} \tag{234}$$

So, substituting for $\gamma$ using (232), now we compute

$$\begin{aligned}
\left\| \left( A - r I_{op} \right) \gamma \right\|^2 &= \lim_{m \to \infty} \left\| \left( A - r I_{op} \right) A\alpha_m \right\|^2 = \lim_{m \to \infty} \left\| AA - rA\alpha_m \right\|^2 \\
&= \lim_{m \to \infty} \left\| A \left( A - r I_{op} \right) \alpha_m \right\|^2 \\
&\leq \|A\|^2 \lim_{m \to \infty} \left\| \left( A - r I_{op} \right) \alpha_m \right\|^2 = \|A\|^2 \lim_{m \to \infty} \left( \left( A - r I_{op} \right) \alpha_m, \left( A - r I_{op} \right) \alpha_m \right) \\
&= \|A\|^2 \lim_{m \to \infty} \left[ \left( A\alpha_m, A\alpha_m \right) + r^2 \left( \alpha_m, \alpha_m \right) - r \left( \alpha_m, A\alpha_m \right) - r \left( A\alpha_m, \alpha_m \right) \right] \\
&= \|A\|^2 \lim_{m \to \infty} \left[ \left\| A\alpha_m \right\|^2 + r^2 \left\| \alpha_m \right\|^2 - 2r \left( \alpha_m, A\alpha_m \right) \right] \\
&\leq \|A\|^2 \lim_{m \to \infty} \left[ \|A\|^2 \left\| \alpha_m \right\|^2 + r^2 \left\| \alpha_m \right\|^2 - 2r \left( \alpha_m, A\alpha_m \right) \right] \\
&= \|A\|^2 \left[ r^2 + r^2 - 2r^2 \right] \\
&= 0
\end{aligned} \tag{235}$$

where we have used $\|A\alpha\| \leq \|A\| \|\alpha\|$ (Eq. (49)), $(\alpha_m, A\alpha_m) = (A\alpha_m, \alpha_m)$ by Hermiticity, $\|A\|^2 = r^2$ from (231), $\lim_{m \to \infty} (\alpha_m, A\alpha_m) = r$ from Eq. (233), and $\|\alpha_m\|^2 = 1$ by definition.

Hence, $\gamma$ is indeed an eigenvector, which proves $r$ is an eigenvalue. Note that this argument proves that one or other or possibly both of $\|A\|$ or $-\|A\|$ is an eigenvalue of $A$, but at least one of them is.

## 11.16.2 Inductive construction of a set of eigenvectors and eigenvalues

With the above procedure, we have found a first eigenvalue $r_1$, with either $r_1 = \|A\|$ or $r_1 = -\|A\|$. (If both of these are eigenvalues, then we choose just one of them for the moment – say the positive one, for definiteness.) $r_1$ may have more than one linearly independent eigenvector associated with it (i.e., it may be a degenerate eigenvalue), but we know that number $m_1$ of linearly independent eigenvectors is finite by (209) above, and $m_1$ is the dimensionality of the space $M_1$ spanned by these eigenvectors. We can therefore construct an orthogonal set of basis functions $\{\beta_1, \ldots, \beta_{m_1}\}$ for the space $M_1$ (and these are also eigenvectors for this eigenvalue). (It is possible, of course, that the eigenvalue $r_1$ is not degenerate, in which case $m_1 = 1$, and there is only one basis function in the set, but for generality we keep the full notation to cover the non-degenerate and degenerate cases.) For notational reasons, the operator $A$ will now also be called $A_1$ (i.e., $A \equiv A_1$).

A subspace $M$ of a Hilbert space is called an *invariant subspace* of a linear operator $B$ if, for any vector $\alpha$ in $M$, $B\alpha$ is also a vector in $M$. We can see that $M_1$ is an invariant subspace of $H$ for the operator $A_1$. Now, when $A_1$ acts on any linear combination of the vectors $\{\beta_1, \ldots, \beta_{m_1}\}$, the result is a vector in the same space $M_1$ since $A_1 \beta_j = r_1 \beta_j$ for any $\beta_j$ in the set. We can then usefully define another operator

$$A_2 = A_1 - \sum_{j=1}^{m_1} r_1 \beta_j \left(\beta_j, \cdot\right) \tag{236}$$

(Here in the mathematical notation $\left(\beta_j, \cdot\right)$ means that, when the operator acts on a vector, we substitute that vector for the dot "$\cdot$" in this inner product expression. This is generally much clearer in Dirac notation, where we would write Eq. (236) as $A_2 = A_1 - \sum_{j=1}^{m_1} r_1 |\beta_j\rangle\langle\beta_j|$, though we will retain the mathematical notation in this proof.)

Now, $A_2$ acting on any vector on $H$, only generates vectors that are orthogonal to all the $\beta_j$, which means that any eigenvectors of $A_2$ are also orthogonal to these $\{\beta_j\}$ and the associated eigenvalues must all be different from $r_1$. So now we can repeat the process we used with $A_1$ to find now a (largest magnitude) eigenvalue $r_2$ of $A_2$ with an associated set of $m_2$ orthogonal vectors $\{\beta_{m_1+1}, \ldots, \beta_{m_1+m_2}\}$. Note that necessarily $|r_2| = \|A_2\| \leq r_1$; it is possible that, if both $+\|A_1\|$ and $-\|A_1\|$ were eigenvalues of $A_1$, $r_2$ is now the "other" one of those, and hence is of equal magnitude to $r_1$. Otherwise, it must be some (real) number of smaller magnitude. Hence, we have now found a second set of eigenfunctions $\{\beta_{m_1+1}, \ldots, \beta_{m_1+m_2}\}$, all orthogonal to the first set $\{\beta_1, \ldots, \beta_{m_1}\}$ and with a different eigenvalue $r_2$. We proceed similarly, with

$$A_{n+1} = A_n - \sum_{j=k_n}^{k_n+m_n} r_n \beta_j \left(\beta_j, \cdot\right) \tag{237}$$

where $k_n$ is the sum of all the preceding $m_q$, i.e.,

$$k_n = \sum_{q=1}^{n-1} m_q \tag{238}$$

with, formally, $k_1 = 0$, and

$$|r_1| \geq |r_2| \geq |r_3| \geq \cdots \tag{239}$$

Note, incidentally, that this means that the eigenvectors and eigenvalues can be found progressively by a variational approach; we would choose a (normalized) test vector $\theta$ constrained to be orthogonal to all preceding eigenvectors, and varying $\theta$ to give the largest possible value of the inner product $(\theta, A\theta)$, and the resulting vector would be the "next" eigenvector, with an associated eigenvalue equal to the resulting maximized value of $(\theta, A\theta)$. Though it might be unlikely that we would in practice use such a variational technique for calculations, this point is conceptually and physically important in establishing eigenvectors as the ones that maximize the inner product $(\theta, A\theta)$. We state this formally above as (114).

If we make a notational change to write the eigenvalues with the same index $j$ as used for the eigenvectors, with the understanding that the eigenvalue $r_j$ is whatever one is associated with the eigenvector $\beta_j$, then we can concatenate all the expressions as in Eq. set to give

$$A_{n+1} = A - \sum_{j=1}^{k_{n+1}} r_j \beta_j (\beta_j, \cdot) \tag{240}$$

Now, we know by the same procedure that

$$\|A_{n+1}\| = |r_{n+1}| \tag{241}$$

Hence, we know that

$$\left\| A - \sum_{j=1}^{k_{n+1}} r_j \beta_j (\beta_j, \cdot) \right\| = |r_{n+1}| \tag{242}$$

and so

$$\lim_{n\to\infty} \left\| A - \sum_{j=1}^{k_{n+1}} r_j \beta_j (\beta_j, \cdot) \right\| = \lim_{n\to\infty} |r_{n+1}| \tag{243}$$

But we know that the eigenvalues of a compact operator must tend to zero as $n \to \infty$ (Eq. (104) (OE5) as proved above in 11.14 "Proof (14) that the eigenvalues of Hermitian operator on an infinite dimensional space tend to zero"), and so we have proved that

$$A = \sum_{j=1}^{\infty} r_j \beta_j (\beta_j, \cdot) \tag{244}$$

where the sum converges in the operator norm. We restate this representation of the operator $A$ above as (112). This also proves that the set of eigenfunctions of a compact operator are complete for describing the effect of the operator on any vector. We state this formally above as (111).

If all the eigenvalues are non-zero, then the set will be complete for the Hilbert space $H$. If not, then we can extend the set by Gram-Schmidt orthogonalization to complete it.

## 11.17 Proof (17) of the equivalence of Dirac and matrix forms of SVD

Formally, we can write a matrix that is diagonal on some basis $\{\gamma_j\}$ as

$$D_{diag} = \sum_j s_j |\gamma_j\rangle\langle\gamma_j| \tag{245}$$

where $s_j$ are the diagonal elements, and we can define two matrices

$$U = \sum_p |\psi_p\rangle\langle\gamma_p| \quad \text{and} \quad V = \sum_q |\phi_q\rangle\langle\gamma_q| \tag{246}$$

It is straightforward to check that these are unitary as in Eq. (140); specifically,

$$U^\dagger U = \left(\sum_p |\psi_p\rangle\langle\gamma_p|\right)^\dagger \left(\sum_q |\psi_q\rangle\langle\gamma_q|\right) = \sum_{p,q} |\gamma_p\rangle\langle\psi_p|\psi_p\rangle\langle\gamma_q|$$
$$= \sum_{p,q} |\gamma_p\rangle\delta_{pq}\langle\gamma_q| = \sum_p |\gamma_p\rangle\langle\gamma_p| = I_{op} \tag{247}$$

and similarly for $V^\dagger V$. Then

$$A = \mathsf{V}\mathsf{D}_{diag}\mathsf{U}^\dagger = \left(\sum_q \left|\phi_q\right\rangle\left\langle\gamma_q\right|\right)\left(\sum_j s_j \left|\gamma_j\right\rangle\left\langle\gamma_j\right|\right)\left(\sum_p \left|\gamma_p\right\rangle\left\langle\psi_p\right|\right)$$

$$= \sum_{q,j,p} \left|\phi_q\right\rangle \delta_{qj} s_j \delta_{jp} \left\langle\psi_p\right| = \sum_j s_j \left|\phi_j\right\rangle\left\langle\psi_j\right| \tag{248}$$

which is identical to Eq. (136), so these two approaches are equivalent.

# 13 Index of definitions